\documentclass[11pt, oneside, letterpaper]{amsart} 

\usepackage[margin=1.15in, left=1.2in, right=1.2in]{geometry}
\usepackage{amsmath,amsfonts,amssymb,amscd,amsthm,amsbsy,epsf,graphics} 
\usepackage{empheq}
\usepackage{stmaryrd}
\usepackage{setspace}
\usepackage{comment}
\usepackage[dvipsnames]{xcolor}
\usepackage[colorlinks,linkcolor=blue,citecolor=blue]{hyperref} 
\usepackage{epsfig} 
\usepackage{subcaption}
\usepackage{tikz-cd}
\usepackage{bm}
\usepackage{graphicx}  
\usepackage{lipsum}
\usepackage{enumitem} % remove indentation of items
\UseRawInputEncoding

\usepackage{todonotes} % use \todo to make marginal comments.
\usepackage[square,numbers]{natbib}

\makeatletter
\def\l@subsection{\@tocline{2}{0pt}{4pc}{5pc}{}}
\let\oldtocsection=\tocsection
\let\oldtocsubsection=\tocsubsection
\let\oldtocsubsubsection=\tocsubsubsection
 
\renewcommand{\tocsection}[2]{\hspace{0em}\oldtocsection{#1}{#2}}
\renewcommand{\tocsubsection}[2]{\hspace{0em}\oldtocsubsection{#1}{#2}}
\renewcommand{\tocsubsubsection}[2]{\hspace{2em}\oldtocsubsubsection{#1}{#2}}

\newcommand{\Abold}{\bm{A}^+}
\newcommand{\ared}[1]{\bm{A}_{{\rm red},\, #1}}

\newtheorem{theorem}{Theorem}[section] 
\newtheorem*{theorem*}{Theorem}
\newtheorem{corollary}[theorem]{Corollary}
\newtheorem*{corollary*}{Corollary}
\newtheorem{lemma}[theorem]{Lemma}
\newtheorem{proposition}[theorem]{Proposition}
\newtheorem*{proposition*}{Proposition}

\newtheorem{question}[theorem]{Question}
\newtheorem*{question*}{Question}

\newtheorem*{problem*}{Problem}

\theoremstyle{definition}

\newtheorem{remark}[theorem]{Remark} 
\newtheorem*{remark*}{Remark} 
\newtheorem{remarks}[theorem]{Remarks}

\newtheorem*{acknowledgement*}{Acknowledgements}

\theoremstyle{plain}
\newtheorem{innercustomthm}{Theorem}
\newenvironment{customthm}[1]
  {\renewcommand\theinnercustomthm{#1}\innercustomthm}
  {\endinnercustomthm}

\newtheorem{innercustomprop}{Proposition}
\newenvironment{customprop}[1]
{\renewcommand\theinnercustomprop{#1}\innercustomprop}
{\endinnercustomprop}

\newtheoremstyle{cases}
  {6pt plus 3 pt}%       Space above
  {2pt}%       Space below
  {\bfseries}   %       Body font
  {}%          Indent amount (empty = no indent, \parindent = para indent)
  {\bfseries}% Thm head font
  {.}%         Punctuation after thm head
  {.5em}%      Space after thm head: " " = normal interword space;
  {}%          Thm head spec (can be left empty, meaning `normal')
\theoremstyle{cases}

\numberwithin{subcase}{case} 
\numberwithin{subsubcase}{subcase}
\numberwithin{equation}{subsection} 
\graphicspath{{figures/}} 
\title[Euler classes of foliations]{On $3$-manifolds admitting co-orientable taut foliations, but none with vanishing Euler class}

\author[Steven Boyer]{Steven Boyer}
\thanks{Steven Boyer was partially supported by NSERC grant RGPIN 2024-06196}
\address{D\'epartement de Math\'ematiques, Universit\'e du Qu\'ebec \`a Montr\'eal, 201 President Kennedy Avenue, Montr\'eal, Qc., Canada H2X 3Y7.}
\email{boyer.steven@uqam.ca}
\urladdr{http://www.cirget.uqam.ca/boyer/boyer.html}

\author[Cameron McA. Gordon]{Cameron McA. Gordon} 
\address{Department of Mathematics, University of Texas at Austin, 1 University Station, Austin, TX 78712, USA.}
\email{gordon@math.utexas.edu}

\author[Ying Hu]{Ying Hu}
\thanks{Ying Hu was partially supported by NSF grant DMS-2409398}
\address{Department of Mathematical and Statistical Sciences, University of Nebraska Omaha, 6001 Dodge Street, Omaha, NE 68182-0243, USA.}
\email{yinghu@unomaha.edu}
\urladdr{https://yinghu-math.github.io}

\author[Duncan McCoy]{Duncan McCoy}
\thanks{Duncan McCoy was partially supported by NSERC grant  RGPIN-2020-05491 and a Canada Research Chair}
\address{D\'epartement de Math\'ematiques, Universit\'e du Qu\'ebec \`a Montr\'eal, 201 President Kennedy Avenue, Montr\'eal, Qc., Canada H2X 3Y7.}
\email{mc\_coy.duncan@uqam.ca}
\urladdr{https://sites.google.com/view/duncanmccoy}



\begin{document}

\begin{abstract}
In this article, we construct infinitely many (small Seifert fibred, hyperbolic and toroidal) rational homology $3$-spheres that admit co-orientable taut foliations, but none with vanishing Euler class. In the context of the $L$-space conjecture, these examples provide rational homology $3$-spheres that admit co-orientable taut foliations (and hence are not $L$-spaces) and have left-orderable fundamental groups, yet none of the left orders arise directly from the universal circle actions associated to co-orientable taut foliations.

The hyperbolic and non-Seifert toroidal examples are obtained from Dehn surgeries on knots in the $3$-sphere and use Heegaard Floer homology to obstruct the existence of a co-orientable foliation with vanishing Euler class. For the Seifert fibred case, we establish necessary and sufficient conditions for the Euler class of the normal bundle of the Seifert fibration to vanish. Moreover, when the base orbifold is hyperbolic, we also provide a second proof of this condition from the viewpoint of discrete faithful representations of Fuchsian groups.
\end{abstract}

\maketitle

\setcounter{tocdepth}{1}
{
\parskip=.2em
\tableofcontents
} 

\section{Introduction}

A long-standing 
problem in $3$-manifold topology is to determine which classes in the second cohomology group of a closed, oriented $3$-manifold $M$ are the Euler classes of tangent bundles of co-oriented taut foliations on $M$. Such classes are constrained by the  fact that every orientable $3$-manifold is parallelizable, which implies that the Euler class of a co-oriented plane field over $M$ is an even class (see \cite[Proposition 2.3]{Hu23} for instance).  If we drop the tautness condition, every even class is realised, since each co-oriented plane field on $M$ is homotopic to an integrable one \cite{Thurston75}.  This no longer holds when one requires that the foliation be taut.  Indeed, only finitely many homotopy classes of co-oriented plane fields can be realised as the tangent plane field of a taut foliation \cite[Corollary 1.5]{KM97}, \cite[Corollary 6.18]{Gab99}. Determining which cohomology classes can occur as Euler classes of co-orientable taut foliations remains open. 

Much work on the problem has been done when $b_1(M)>0$. In this case, Gabai has shown that taut foliations always exist \cite{Gab83}. Further, results of Roussarie and Thurston show that the dual Thurston norm of the Euler class of any co-oriented taut foliation is at most one \cite{Roussarie, Thurston86}. Thurston also conjectured that any class in $H^2(M)$ of dual Thurston norm one can be realised as the Euler class of a co-oriented taut foliation. However, this conjecture is now known to be false, following counterexamples constructed by Yazdi and Liu \cite{Yazdi20, Liu24}. At the other extreme, one may ask whether the zero class can be realised as the Euler class of a co-oriented taut foliation. See \cite[Question 9.1 and Question 9.2]{Yazdi20} and the related discussions therein. Also see \cite{FLY25}.

In this article, we consider the case $b_1(M) = 0$; that is, $M$ is a rational homology sphere or equivalently, $H^2(M)$ is finite. This setting differs from the case $b_1(M) > 0$ in two fundamental ways. First, the existence of co-orientable taut foliations is no longer guaranteed. Second, the dual Thurston norm vanishes, since $H^2(M)$ is finite, so the usual techniques used to analyze the Euler class when $b_1(M) > 0$ do not apply. Nevertheless, the zero class remains a distinguished element of $H^2(M)$. This leads naturally to the following question (cf. \cite[Question 9.2]{Yazdi20}):

\begin{question}
\label{que: basic qn}
If a closed, orientable rational homology sphere admits a co-orientable taut foliation, does it admit one whose Euler class is zero? \footnote{Given a co-orientable taut foliation, the Euler class of its tangent bundle is well-defined up to sign. So it makes sense to discuss if a co-orientable taut foliation has vanishing Euler class.}
\end{question}

Question \ref{que: basic qn} takes on further significance in the context of the $L$-space conjecture \cite{BGW13, Juhasz2015}, which predicts that if a closed, orientable, prime $3$-manifold $M$ admits a co-orientable taut foliation, then its fundamental group is left-orderable. The conjecture is known to hold when $b_1(M)>0$ \cite{Gab83, BRW05}. On the other hand, when $M$ is a rational homology sphere, one of the main tools for proving this is to use Thurston's universal circle action, which to each co-oriented taut foliation $\mathcal{F}$ on an oriented rational homology $3$-sphere $M$ associates a non-trivial representation 
$$\rho_\mathcal{F}: \pi_1(M) \to \mbox{Homeo}_+(S^1).$$
See \cite{CD03} for the details. If one can show that this representation lifts to the universal covering group $\mbox{Homeo}_{\mathbb Z}(\mathbb R)$ of $\mbox{Homeo}_+(S^1)$, then it follows from \cite[Theorem 1.1]{BRW05} that the fundamental group of $M$ is left-orderable. However, such a lift exists if and only if the Euler class of $\rho_\mathcal{F}$ vanishes. Since this Euler class coincides with the Euler class of the tangent bundle of $\mathcal{F}$ \cite{BH19}, which we denote by $e(\mathcal{F})$, we are naturally led to Question \ref{que: basic qn} 

For the reasons discussed above, Question \ref{que: basic qn} has been raised often over the years. Its answer is yes when $H^2(M)$ is a (possibly trivial) group consisting of 2-torsion, since the only even class is zero in this case. And though it would be surprising if the answer were yes in general, no examples of manifolds for which the answer is no have appeared in the literature. In this article, we show that such examples are abundant, and ubiquitous amongst manifolds obtained by surgery on knots in the $3$-sphere (Theorem \ref{thm:hyperbolic_examples}) and amongst Seifert fibred rational homology $3$-spheres (Theorem \ref{thm:all_homology_groups}). Further, many of our examples have left-orderable fundamental groups, even though no left-order on the group arises from a lift of a universal circle representation. Thus we ask the following question.

\begin{question}
Is there a general procedure for constructing, from a closed, orientable, prime 3-manifold which admits a co-orientable taut foliation, a non-trivial action of its fundamental group on the reals by orientation-preserving homeomorphisms?
\end{question}

\subsection{Surgery on knots in the $3$-sphere} 
\label{subsec: intro surgery on knots}
The examples derived from surgery are obtained by applying recent work of Lin, who showed that if a $3$-manifold $M$ admits a co-oriented taut foliation $\mathcal{F}$, then there is a $\mathrm{spin}^c$-structure $\mathfrak{s}$ on $M$ whose associated plane field is $T\mathcal{F}$ and a direct $\mathbb{F}$ summand of the $\mathbb{F}[U]$-module $HF_{\mathrm{red}}(M,\mathfrak{s})$ \cite{Lin24}.

{\bf Convention}: When referring to $p/q$-surgery on a knot, we assume that $p$ and $q$ are coprime integers with $p > 0$.  We denote the manifold obtained by $p/q$-surgery on a knot $K$ by $K(p/q)$.

\begin{theorem}
\label{thm: p/q-surgery}
Let $K$ be a knot in the $3$-sphere. If $q$ is even and $|p/q| > 2g(K) - 1$, then $K(p/q)$ does not admit a co-orientable taut foliation with Euler class zero.
\end{theorem}

 In \cite[Theorem 1.12(1)]{Hu23}, it was shown that if $|p/q| > 2g(K)$ and $p/q$ is not an integer, then any co-orientable taut foliation on $K(p/q)$ that is transverse to the core of the surgery solid torus cannot have vanishing Euler class. The above theorem removes the transversality condition under the assumption that $q$ is even. 
 
Thus by Theorem \ref{thm: p/q-surgery}, the existence of a knot for which some large surgery with even denominator admits a co-orientable taut foliation shows that the answer to Question \ref{que: basic qn} is no.  We give some examples in the corollary below.

\begin{corollary}
\label{cor: fdtc 0}
Suppose that $K$ is a non-trivial knot in the $3$-sphere which is either 
\vspace{-.2cm}
\begin{enumerate}
    \setlength{\itemsep}{0.2em}
\item[{\rm (1)}] fibred with zero fractional Dehn twist coefficient, or
\item[{\rm (2)}] alternating but not a $(2,q)$-torus knot, or

\item[{\rm (3)}] a Montesinos knot but not a $(2,q)$-torus knot or a $(-2,3,2n+1)$-pretzel knot.

\end{enumerate}
Then for all slopes $p/q$ with $q$ even and $|p/q| > 2g(K) -1$, $K(p/q)$ admits a co-orientable taut foliation, but none with zero Euler class. 
\end{corollary}

Part (1) of Corollary \ref{cor: fdtc 0} follows from Theorem \ref{thm: p/q-surgery} together with Roberts' result \cite{RobertsSurfacebundle2}, while parts (2) and (3) follow from  Theorem \ref{thm: p/q-surgery} and work of Delman and Roberts.

As we mentioned above, examples of manifolds having left-orderable fundamental groups and admitting co-orientable taut foliations, but none with vanishing Euler class, can also be constructed, thus showing that no left-orders on the group arise from lifts of universal circle representations. Here are the simplest examples. More examples are provided in Theorem \ref{thm:hyperbolic_examples} and Theorem \ref{thm:all_homology_groups} below.
 
\begin{itemize}

\item  The figure eight knot $K$ is a fibred hyperbolic knot of genus one with zero fractional Dehn twist coefficient and each rational surgery on it has a left-orderable fundamental group \cite{Zung24}. Hence, Corollary \ref{cor: fdtc 0} implies that for all slopes $p/q$ with $q$ even and $|p/q| > 1$, $K(p/q)$ has a left-orderable fundamental group and admits a co-orientable taut foliation, but no co-orientable taut foliation with zero Euler class. 

\vspace{.2cm} \item Suppose that $K$ is the right-handed trefoil knot. Then $K(p/q)$ admits a co-orientable taut foliation (equivalently, has a left-orderable fundamental group \cite{BRW05}) if and only if $p/q < 1$ \cite{JN85b, Na94}. It follows from Lemma \ref{lem: same Euler class} and \cite[Theorem~1.7]{Hu23} that $K(p/q)$ admits a co-orientable taut foliation with zero Euler class if and only if $p/q < 1$ and $|q| \equiv 1 \pmod{p}$.

\end{itemize}
All the examples in the first family are hyperbolic, while those in the second are Seifert fibred. We can also use Theorem \ref{thm: p/q-surgery} to produce toroidal examples which are not Seifert fibred. 

Define the {\it JSJ graph} of a knot $K$ to be the rooted tree dual to the JSJ tori of the exterior $X(K)$ of $K$, where the root vertex corresponds to the piece of the JSJ decomposition containing $\partial X(K)$. Any finite rooted tree can arise as the JSJ graph of a knot. We say that a JSJ graph is a {\it rooted interval} when it is homeomorphic to an interval with root at an endpoint. Since non-integral surgeries on knots in $S^3$ are irreducible \cite{GLu87}, the following corollary follows immediately from Theorem \ref{thm: p/q-surgery} and Theorems 3.13 and 4.3 of \cite{BGH25}. 

\begin{corollary}
\label{cor: non-rooted interval satellite}
Suppose that $K$ is a fibred satellite knot whose $JSJ$ graph is not a rooted interval. Then for all slopes $p/q$ with $q$ even and $|p/q| > 2g(K) - 1$, $K(p/q)$ is toroidal, has a left-orderable fundamental group, and admits a co-orientable taut foliation, but none with zero Euler class. 
\end{corollary}

\begin{remark}
Fibred composite knots satisfy the corollary's hypothesis. Taking $K$ to be a connected sum of torus knots yields examples which are graph manifolds. 
\end{remark}

If $K$ is an $L$-space knot then for slopes $p/q \ge 2g(K)-1$ the manifold $K(p/q)$ is an $L$-space and hence does not admit a co-orientable taut foliation. On the other hand, since $L$-space knots are fibred, Roberts' work \cite{RobertsSurfacebundle2} shows that $K(p/q)$ has a co-orientable taut foliation if $p/q < 1$. We establish conditions under which negative surgeries on an $L$-space knot cannot admit a co-orientable taut foliation with Euler class zero. These conditions are phrased in terms of congruence relations involving the torsion coefficients of the Alexander polynomial (Theorem \ref{thm:neg_L_space_obstruction} and Corollary \ref{cor: q arbitrary t0}), as well as the bridge index (Corollary \ref{cor:bridge_bound}).

By performing negative surgeries on the $(-2, 3, 2n+1)$-pretzel knots, Theorem~\ref{thm:neg_L_space_obstruction} provides the following examples.

 \begin{theorem}
 \label{thm:hyperbolic_examples}
     For any integer $p\geq 7$ there are infinitely many hyperbolic 3-manifolds $M$ such that $H^2(M)\cong \mathbb{Z}/p$ and $M$ admits a co-orientable taut foliation, but none with Euler class zero. Further, $\pi_1(M)$ is left-orderable.
 \end{theorem}
 
By way of contrast, note that Theorem~\ref{thm: p/q-surgery} cannot be extended to large surgeries on all knots without some hypotheses on the slope $p/q$. For instance, if $M$ is a rational homology sphere for which there is a unique $\mathrm{spin}^c$-structure $\mathfrak{s}_0$ such that $HF_{\mathrm{red}}(M,\mathfrak{s}_0)$ is non-zero, then any co-orientable taut foliation on $M$ has Euler class zero (see Proposition \ref{prop:single_spinc_structure}). As a consequence, we deduce,

\begin{customprop}{\ref{prop:genus_one_surgeries}}
Let $M$ be a rational homology sphere obtained by $p$-surgery on a knot $K$ is the $3$-sphere, for some $p \in  \mathbb Z$. If 
$K$ has genus one, or $K$ is an almost $L$-space knot \footnote{$K$ is an almost $L$-space knot if it has a \emph{positive} almost $L$-space surgery of slope $p/q \geq 2g(K)-1$ \cite{BS24}.} and $p\geq 2g(K)-1$, then any co-orientable taut foliation on $M$ has Euler class zero. 
\end{customprop}

Finally, recalling the theorem of Hu from \cite{Hu23} mentioned above (also see Corollary~\ref{cor: q arbitrary t0} and Remark \ref{rem: t_0 condition}), we ask the following question.

\begin{question}
\label{que: q=1}
In Theorem \ref{thm: p/q-surgery}, can the assumption that $q$ be even be replaced by the assumption $|q| > 1$?
\end{question}

\subsection{Seifert fibred manifolds}
We denote a closed oriented Seifert fibred manifold $M$ with orientable base orbifold by $M\bigl(g;\, b, \tfrac{a_1}{b_1}, \ldots, \tfrac{a_n}{b_n}\bigr)$, where $b \in \mathbb{Z}$ and $b_i \geq 2$ for $i = 1, \ldots, n$. See \S\ref{subsection: notation} for a detailed explanation of the notation. We denote the base orbifold by $\mathcal{B}$; this is a closed orientable surface $|\mathcal{B}|$ of genus $g \ge 0$ with $n \ge 0$ cone points of orders $b_1,...,b_n$. Recall that the Euler characteristic $\chi(\mathcal{B})$ of $\mathcal{B}$ is $\chi(|\mathcal{B}|) - n + \sum_{i=1}^n \frac{1}{b_i}$, and the Euler number of $M$ is $e(M) = b + \sum_{i=1}^n \frac{a_i}{b_i}$.

A Seifert fibred rational homology $3$-sphere $M$ which admits a co-orientable taut foliation $\mathcal{F}$ is necessarily of the form $M\bigl(0;\, b, \tfrac{a_1}{b_1}, \ldots, \tfrac{a_n}{b_n}\bigr)$ \cite[Proposition 6.1]{BC17}. Further, $\mathcal{F}$ is necessarily {\it horizontal}, that is, transverse to the Seifert fibres of $M$ (see Lemma \ref{lem: same Euler class}). As such, the tangent bundle of $\mathcal{F}$ is isomorphic to the normal bundle $\nu_M$ to the Seifert fibres of $M$. In \S \ref{sec: Euler class normal bundle} we determine necessary and sufficient conditions for the Euler class of $\nu_M$ to vanish (Theorem \ref{thm: Euler class zero}), which leads to the following theorem.

\begin{theorem}
\label{thm: main application}
Suppose that $M =  M\bigl(0;\,b, \tfrac{a_1}{b_1}, \ldots, \tfrac{a_n}{b_n}\bigr)$ is a rational homology $3$-sphere which admits a co-orientable taut foliation $\mathcal{F}$. Then $\mathcal{F}$ has zero Euler class if and only if there exists $m \in \mathbb{Z}$ such that
\vspace{-.2cm}
\begin{enumerate}

\vspace{.2cm} \item[\rm (1)] $m a_i \equiv  1 \, (\text{{\rm mod} }b_i), \text{ for } i =  1, \cdots, n, \; \text{and}$

\vspace{.2cm} \item[\rm (2)] $ me(M) = \chi(\mathcal{B})$.

\end{enumerate}

\end{theorem}

\begin{corollary}
    \label{cor: necessary condition intro}
Given a closed oriented Seifert fibred manifold $M$ which admits a co-orientable taut foliation whose Euler class vanishes, then one of the following holds:
\vspace{-.2cm}
\begin{enumerate}
\vspace{.2cm} \item[\rm (1)]   $M$ is an orientable $S^1$-bundle over the torus, 
\vspace{.2cm} \item[\rm (2)]  $e(M) = \chi(\mathcal{B}) = 0$, 
\vspace{.2cm} \item[\rm (3)] $ e(M)\neq 0, \, \chi(\mathcal{B}) \neq 0, $ and $$\frac{\chi(\mathcal{B})}{e(M)}\in \mathbb{Z}.$$ 

\end{enumerate}
 
\end{corollary}

By abelianizing the presentation (\ref{eqn: pres 2}) of $\pi_1(M)$ obtained by viewing $M$ as doing a Dehn filling of a trivial $S^1$-bundle over an $n$-punctured orientable surface of genus $g$ as in \S \ref{subsection: notation}, one can show that if $e(M) \ne 0$ and $T_1(M)$ denotes the torsion subgroup of $H_1(M)$, then 
$$|T_1(M)| = b_1 \cdots b_n | e(M)|,$$
where we take the convention that $b_1b_2 \cdots b_n = 1$ when $n = 0$. Thus, it follows from Theorem \ref{thm: main application} that:  

\begin{corollary}
\label{cor: 2}
Suppose that $M =  M\bigl(0;\,b, \tfrac{a_1}{b_1}, \ldots, \tfrac{a_n}{b_n}\bigr)$ is a rational homology $3$-sphere which admits a co-orientable taut foliation $\mathcal{F}$ with zero Euler class. If $e(M) \ne 0$ then the integer $b_1\ldots b_n \chi(\mathcal{B})$ is an integer multiple of $|H_1(M)|$. 
\end{corollary}

When the base orbifold $M = M\bigl(g;\,b, \tfrac{a_1}{b_1}, \ldots, \tfrac{a_n}{b_n}\bigr)$ is hyperbolic, we provide a second proof of Theorem \ref{thm: main application} from the viewpoint of discrete faithful representations of Fuchsian groups. Here is an interesting by-product of the analysis.

\begin{theorem}
\label{thm: e = 0 means ctf}
Let $M$ be a closed oriented Seifert fibred manifold with orientable base orbifold $\mathcal{B}$. If $\mathcal{B}$ is hyperbolic or Euclidean, and $e(\nu_M) = 0$, then $M$ admits a horizontal foliation.
\end{theorem}
 
As an application of Theorem \ref{thm: main application} we show that such examples are actually abundant in the presence of higher order torsion.  

\begin{theorem}\label{thm:all_homology_groups}
    Let $G$ be a finite abelian group containing an element of order at least three. Then there exist infinitely many Seifert fibred manifolds $M$ such that $H_1(M)\cong G$ and $M$ admits a co-orientable taut foliation, but none with Euler class zero. Further, $\pi_1(M)$ is left orderable.
\end{theorem}

\begin{remark}
It is known that Seifert fibred manifolds satisfy the $L$-space conjecture \cite{BRW05, LS07, BGW13}, so all such manifolds which admit co-orientable taut foliations have left-orderable fundamental groups even though it is often the case that this cannot be shown using universal circle representations. One uses, rather, the fact that co-orientable taut foliations on Seifert fibred rational homology spheres are $\mathbb R$-covered. 
\end{remark}

We will present more examples of Seifert fibred manifolds that admit co-orientable taut foliations, but none with vanishing Euler class, in 
\S\ref{subsec: hyperbolic}. In fact, our computations in \S\ref{sec: horizontal foliation} suggest that most closed 
Seifert fibred manifolds have nontrivial normal bundle, including the case when the base orbifold is Euclidean (i.e.\ $\chi(\mathcal{B}) = 0$).  This may seem counterintuitive, since we have pointed out that the normal bundle is always trivial when the base $\mathcal{B}$ is the torus. In contrast, closed Seifert fibred rational homology spheres with spherical base orbifolds have finite fundamental groups and hence cannot admit taut foliations. Yet, in \S\ref{subsec: spherical}, we will construct examples of such manifolds for which 
$e(\nu_M)$ vanishes.

\subsection*{Plan of the paper}
In \S \ref{sec: surgery examples}, we apply Heegaard Floer homology to produce hyperbolic $3$-manifolds obtained by Dehn filling on knots in the $3$-sphere that admit co-orientable taut foliations, but none with vanishing Euler class. All results stated in \S \ref{subsec: intro surgery on knots} are proved in this section.   In \S\ref{sec: Euler class normal bundle} we explain our notation for Seifert manifolds and prove Theorem~\ref{thm: Euler class zero}, which characterizes those Seifert manifolds for which normal bundle $\nu_M$ to the Seifert fibres is trivial. Theorem \ref{thm: main application} then follows from this and Lemma \ref{lem: same Euler class}. In \S\ref{sec: normal} we study the Euler class of the normal bundle $\nu_M$ through the use of discrete faithful $PSL(2,\mathbb{R})$-representations of $\pi_1(\mathcal{B})$ in the case that the base orbifold $\mathcal{B}$ is hyperbolic. As a consequence of this approach, we provide a different proof of Theorem \ref{thm: Euler class zero} when the base is hyperbolic and prove the hyperbolic case of Theorem \ref{thm: e = 0 means ctf}. Finally, in \S\ref{sec: horizontal foliation} we provide examples of Seifert fibred 
manifolds whose normal bundles are (non-)trivial. We classify which occurs when $\mathcal{B}$ is spherical in \S \ref{subsec: spherical} and is Euclidean in \S \ref{subsec euclidean}. The Euclidean case of Theorem \ref{thm: e = 0 means ctf} is dealt with at the end of \S \ref{subsec euclidean}. The case that $\mathcal{B}$ is hyperbolic is analysed in \S \ref{subsec: hyperbolic}.

\subsection*{Acknowledgements}
We thank Nathan Dunfield for bringing up the question of finding examples of rational homology spheres that admit co-oriented taut foliations but none with vanishing Euler class, during his lecture at the workshop ``Low Dimensional Topology and Floer Theory'' in Montreal, August 18-29, 2025, which inspired us to investigate the Euler class of horizontal foliations in Seifert fibred manifolds. The second and third named authors would like to thank the Simons Foundation, the CRM, and CIRGET for their support during the preparation of this article.

\section{Examples from surgeries on knots in the $3$-sphere}
\label{sec: surgery examples}
In this section we prove our results on surgeries on knots in the $3$-sphere. 

\subsection{$\mathrm{Spin}^c$ structures and co-orientable taut foliations} 

The set $\mathrm{Spin}^c(M)$  of $\mathrm{spin}^c$ structures $\mathfrak{s}$ on a rational homology $3$-sphere $M$ corresponds to certain classes of pairs $(\xi, \mathfrak{s}_0)$, where $\xi$ is an oriented $2$-plane sub-bundle of the tangent bundle of $M$ and $\mathfrak{s}_0$ is a spin structure on $M$ \cite[proof of Proposition 1]{Gom97}. There is a faithful, transitive action of $H^2(M)$ on $\mathrm{Spin}^c(M)$ \cite[Proposition 1]{Gom97}. The {\it Chern class} of a $\mathrm{spin}^c$ structure $\mathfrak{s}$, denoted $c_1(\mathfrak{s})$, is the Euler class of its oriented plane bundle. 

The set $\mathrm{Spin}(M)$ of $\mathrm{spin}$-structures on $M$ is identified with $\{\xi_0\} \times \mathrm{Spin}(M) \subset \mathrm{Spin}^c(M)$, where $\xi_0$ is a trivial bundle. This is precisely the set of $\mathrm{spin}^c$ structures $\mathfrak{s}$ satisfying $c_1(\mathfrak{s}) = 0$ \cite[Remark (c) page 49]{Gom97}. Reversing the orientation of $\xi$ determines an involution on $\mathrm{Spin}^c(M)$, called {\it conjugation}, whose fixed point set is $\mathrm{Spin}(M)$ \cite[Remark (c) page 49]{Gom97}. Thus there are precisely $|\mathrm{Spin}(M)| = |H_1(M;\mathbb{Z}/2)|$ self-conjugate spin$^c$-structures on $M$. We record these remarks in a lemma. 

 \begin{lemma} 
 \label{lemma: spinc}
      Let $M$ be a rational homology $3$-sphere. Then there are $|H_1(M)|$ spin$^c$-structures on $M$ and $|H_1(M;\mathbb{Z}/2)|$ self-conjugate spin$^c$-structures. Further, a spin$^c$-structure $\mathfrak{s}$ is self-conjugate if and only if $c_1(\mathfrak{s})=0$.  
\end{lemma}

In what follows, we take $\mathbb F$ to be $\mathbb Z/2$. 

\begin{lemma}\label{lem:lin_obstruction}
Let $M$ be a rational homology $3$-sphere admitting a co-oriented taut foliation $\mathcal{F}$. Then $M$ has a $\mathrm{spin}^c$ structure $\mathfrak{s}$ such that $HF_{\mathrm{red}}(M,\mathfrak{s})$ contains a direct $\mathbb{F}$ summand and $c_1(\mathfrak{s})=e(\mathcal{F})$.
\end{lemma}

\begin{proof}
Let $M$ be a rational homology $3$-sphere that admits a co-oriented taut foliation $\mathcal{F}$. The tangent plane field of $\mathcal{F}$ can be perturbed to obtain a tight contact structure $\xi$ whose contact invariant $c(\xi) \neq 0 \in HF_{\mathrm{red}}(M,\mathfrak{s})$, with $c_1(\mathfrak{s}) = e(\mathcal{F})$ \cite{OS04}. Moreover, Lin deduces that the $\mathbb{F}[U]$-module  $HF_{\mathrm{red}}(M,\mathfrak{s})$ contains direct $\mathbb{F}$ summand \cite{Lin24}.   
\end{proof}

Combined, these two lemmas yield the following obstruction to the existence of taut foliations of Euler class zero and non-zero.
\begin{proposition}\label{prop:euler_class_obstruction}
    Let $M$ be a rational homology sphere with a co-oriented taut foliation $\mathcal{F}$.
    \begin{enumerate}
     \setlength\itemsep{0.5em}
        \item[\rm (1)] If $e(\mathcal{F})=0$, then there is a self-conjugate $\mathrm{spin}^c$ structure $\mathfrak{s}$ such that $HF_{\mathrm{red}}(M,\mathfrak{s})$ contains a direct $\mathbb{F}$ summand.
        \item[\rm (2)] If $e(\mathcal{F})\neq 0$, then there is $\mathrm{spin}^c$ structure $\mathfrak{s}$ such that $HF_{\mathrm{red}}(M,\mathfrak{s})$ contains a direct $\mathbb{F}$ summand and $\mathfrak{s}$ is not self-conjugate.
    \end{enumerate}
\end{proposition}
In particular, Heegaard Floer homology can be used to show that certain 3-manifolds, including almost $L$-spaces, do not admit co-orientable taut foliations with non-zero Euler class.
\begin{proposition}\label{prop:single_spinc_structure}
    Let $M$ be a rational homology sphere for which there is a unique $\mathrm{spin}^c$ structure $\mathfrak{s}_0$ for which $HF_{\mathrm{red}}(M,\mathfrak{s}_0)$ is non-zero. Then any co-orientable taut foliation on $M$ has Euler class zero.
\end{proposition}
\begin{proof}
    Since conjugation of spin$^c$-structures induces an isomorphism on Heegaard Floer homology, we see that $\mathfrak{s}_0$ is necessarily self-conjugate. Thus Proposition~\ref{prop:euler_class_obstruction} implies that $M$ cannot admit a co-orientable taut foliation with non-zero Euler class.
\end{proof}

\subsection{Surgery on knots}
In order to apply Proposition~\ref{prop:euler_class_obstruction}, we need examples of manifolds for which we can understand the Heegaard Floer homology. For surgery on knots $K$ in the $3$-sphere, the Heegaard Floer homology can be calculated using the mapping cone formula \cite{OS08-integer, OS11-rational}.

Let $CFK^\infty(K) = \bigoplus_{i,j \in \mathbb{Z}} C\{(i,j)\}$ be the knot Floer complex of $K$.
We will use $\ared{k}(K)$ to denotes the reduced part of
\[
\Abold_k = H_*(C\{i \ge 0 \text{ or } j \ge k\}).
\]

The mapping cone formula expresses the reduced Heegaard Floer homology of surgeries on $K$ in terms of the $\ared{k}(K)$ \cite{OS08-integer, OS11-rational}. For our purposes, the following calculation will be sufficient.
\begin{lemma}\label{lem:HFred_calculation}
    Let $M$ be obtained by $p/q$-surgery on a knot $K$ of genus $g$ in $S^3$ with $p/q\geq  2g-1$. Then there is an affine map
    \[
\phi : \mathbb{Z}/p \to \mathrm{Spin}^c(M),
\]
such that for $i$ in the range $\frac{q-p}{2}\leq i\leq  \frac{p+q-1}{2}$
\[
HF_{\mathrm{red}}(M, \phi(i)) \cong
\begin{cases}
\ared{\left\lfloor \frac{i}{q} \right\rfloor}(K), & q(1-g) \leq i \leq  qg - 1, \\
0, & \text{otherwise}.
\end{cases}
\]
\end{lemma}
\begin{proof}
 Since $p/q \geq 2g - 1$, the mapping cone formula for computing $HF^+(M)$ in terms of $CFK^\infty(K)$ implies that there is a labelling
\[
\phi : \mathbb{Z}/p \to \mathrm{Spin}^c(M),
\]
such that \cite[Proposition 3.5]{NZ14}
\[
HF_{\mathrm{red}}(M, \phi(i)) \cong \bigoplus_{k \in \mathbb{Z}} \ared{\lfloor \frac{pk+i}{q} \rfloor}(K).
\]
Each value of $i$ (mod $p$) is represented for $i$ in the range $\frac{q-p}{2}\leq i\leq  \frac{p+q-1}{2}$, so as $\ared{j}(K) = 0$ for $|j| \geq  g$ \cite[Lemma 3.2]{NZ14}, we have 
\[
HF_{\mathrm{red}}(M, \phi(i)) \cong
\begin{cases}
\ared{\lfloor \frac{i}{q} \rfloor}(K), & q(1-g) \leq i \leq qg-1, \\
0 & \text{otherwise}.
\end{cases}
\]
\end{proof}
Combined, Proposition~\ref{prop:euler_class_obstruction} and Lemma~\ref{lem:HFred_calculation} quickly yield Theorem~\ref{thm: p/q-surgery}.

\begin{proof}[Proof of Theorem~\ref{thm: p/q-surgery}]
Since $p$ is odd, there is an odd number of spin$^c$-structures on $M$ of which exactly one, call it $\mathfrak{s}_0$, is self-conjugate (Lemma \ref{lemma: spinc}). 
After possibly reversing orientation on $M$ and replacing $K$ by its mirror image, we may assume that $p/q\geq 2g(K)-1$. By Lemma~\ref{lem:HFred_calculation}, we see that for each $j$ such that $|j|<g(K)$, there are precisely $q$ spin$^c$-structures for which $HF_\mathrm{red}(M,\mathfrak{s})\cong \ared{j}(K)$. Thus the number of spin$^c$-structures for which $HF_\mathrm{red}(M,\mathfrak{s})\neq 0$ is divisible by $q$ and hence even. Consequently, we see that there are an odd number of spin$^c$-structures such that $HF_\mathrm{red}(M,\mathfrak{s})=0$. Since conjugation of $\mathrm{spin}^c$ structures induces an isomorphism of Heegaard Floer homology, this implies that $\mathfrak{s}_0$ on $M$ satisfies $HF_\mathrm{red}(M,\mathfrak{s}_0)= 0$. By Proposition~\ref{prop:euler_class_obstruction}, this implies that no co-orientable taut foliation on $M$ can have Euler class zero.
\end{proof}

In the opposite direction, we note that Theorem~\ref{thm: p/q-surgery} cannot be extended to all large surgeries on a knot $K$ without some hypotheses on the slope $p/q$.

\begin{proposition}\label{prop:genus_one_surgeries}
    Let $M$ be a rational homology sphere obtained by $p$-surgery on a knot $K$ is the $3$-sphere, for some $p \in  \mathbb Z$. If 
$K$ has genus one, or $K$ is an almost $L$-space knot and $p\geq 2g(K)-1$, then any co-orientable taut foliation on $M$ has Euler class zero. 
\end{proposition}

\begin{proof}
First suppose that $K$ is an almost $L$-space knot. Then $K(p)$ is an almost $L$-space for each integer $p \geq 2g(K) - 1$ \cite{BS24} and as such, there is a unique $\mathrm{spin}^c$ structure $\mathfrak{s}_0$ for which $HF_{\mathrm{red}}(M,\mathfrak{s}_0) \neq 0$. Hence, Proposition~\ref{prop:single_spinc_structure}  implies that every co-orientable taut foliation on $M$ has Euler class zero.

Next suppose that $K$ is a genus one knot and let $M$ be obtained by $p$-surgery on $K$. Without loss of generality we may assume that $p>0$. In particular, we have $p\geq 2g(K)-1=1$ and we may apply Lemma~\ref{lem:HFred_calculation} to calculate $HF_{\mathrm{red}}(M)$. This shows that $HF_{\mathrm{red}}(M,\phi(i))\cong 0$ for all $i\not\equiv 0 \; (\bmod \; p)$. Thus there is at most one spin$^c$-structure $\mathfrak{s}_0$ with $HF_{\mathrm{red}}(M,\mathfrak{s}_0)\neq 0$. By Proposition~\ref{prop:single_spinc_structure}, this implies that every co-orientable taut foliation on $M$ has Euler class zero. 
\end{proof}
    
\subsection{Negative surgeries on $L$-space knots}
\label{subsec: neg surg lsk}
Let $K$ be a knot in the 3-sphere. Its symmetrized Alexander polynomial can be written in the form
\[
\Delta_K(t)=\sum_{k} a_k t^{k},
\]
where $a_k=a_{-k}$ for all $k$. For $i\geq 0$ we will use $t_i(K)$ to denote $i$th torsion coefficient
\[
t_i(K)=\sum_{k\geq 1}ka_{i+k}.
\]

In what follows we will use $\mathcal{T}(n)$ to denote the cyclic $\mathbb{F}[U]$-module
\[
\mathcal{T}(n)=\frac{\mathbb{F}[U]}{(U^n)}.
\]
\begin{lemma}\label{lem:negative_surgery_calc}
Let $M$ be obtained by $-p/q$-surgery on a non-trivial $L$-space knot $K$ for $p/q>0$. If $q$ is odd, then $M$ has a self-conjugate spin$^c$-structure $\mathfrak{s}_0$ such that
\[
HF_\mathrm{red}(M,\mathfrak{s}_0)\cong \mathcal{T}(t_0) \oplus \bigoplus_{k\geq 1} \mathcal{T}(t_{\lfloor \frac{2pk+q-1}{2q}\rfloor})\oplus \mathcal{T}(t_{\lfloor \frac{2pk+q-1}{2q}\rfloor}).
\]
If $p+q$ is odd, then $M$ has a self-conjugate spin$^c$-structure $\mathfrak{s}_1$ such that
\[
HF_\mathrm{red}(M,\mathfrak{s}_1)\cong \bigoplus_{k\geq 0} \mathcal{T}(t_{\lfloor \frac{2pk+p+q-1}{2q}\rfloor})\oplus \mathcal{T}(t_{\lfloor \frac{2pk+p+q-1}{2q}\rfloor})
\]
\end{lemma}

\begin{proof}
When applied to negative surgeries on $L$-space knots \cite[Lemma~18]{Gainullin}, the mapping cone formula implies that there is a labelling of spin$^c$-structures
\[
\phi: \mathbb{Z}/p \rightarrow \mathrm{Spin}^c(M)
\]
such that
\begin{equation}\label{eq:negative_surgery_formula}
HF_\mathrm{red}(M,\phi(i))\cong \bigoplus_{k\in \mathbb{Z}} \mathcal{T}(t_{|\lfloor \frac{pk+i}{q}\rfloor|}),
\end{equation}
where we are using that the fact that for $L$-space knots, the quantities $H_k$ and $V_k$ are determined by the torsion coefficients of the Alexander polynomial:
\[
V_k=H_{-k}=t_{k}(K).
\]
Using the identity
\[
\left\lfloor \frac{pk+i}{q} \right\rfloor =-\left\lfloor \frac{q-1-i-pk}{q} \right\rfloor  
\]
alongside \eqref{eq:negative_surgery_formula} shows that $HF_\mathrm{red}(M,\phi(i)) \cong HF_\mathrm{red}(M,\phi(\tau(i)))$, where $\tau$ is the involution $\tau(i)=q-1-i\; (\bmod \; p)$ on $\mathbb{Z}/p$. The fixed points of $\tau$ are
\[
\left\{\frac{q-1}{2}, \frac{p+q-1}{2}\right\}\cap \mathbb{Z}.
\]

If $q$ is odd, then we take $\mathfrak{s}_0=\phi(\frac{q-1}{2})$. If $p+q$ is odd, then we take $\mathfrak{s}_1=\phi(\frac{p+q-1}{2})$. Note that $\mathfrak{s}_0$ and $\mathfrak{s}_1$ are both defined if and only if $p$ is even and that precisely one of them is defined if $p$ is odd. For these spin$^c$-structures \eqref{eq:negative_surgery_formula} reduces to 
\[
HF_\mathrm{red}(M,\mathfrak{s}_0)\cong \mathcal{T}(t_0) \oplus \bigoplus_{k\geq 1} \mathcal{T}(t_{\lfloor \frac{2pk+q-1}{2q}\rfloor})\oplus \mathcal{T}(t_{\lfloor \frac{2pk+q-1}{2q}\rfloor})
\]
and
\[
HF_\mathrm{red}(M,\mathfrak{s}_1)\cong \bigoplus_{k\geq 0} \mathcal{T}(t_{\lfloor \frac{2pk+p+q-1}{2q}\rfloor})\oplus \mathcal{T}(t_{\lfloor \frac{2pk+p+q-1}{2q}\rfloor}).
\]
Since $K$ is non-trivial, $t_0(K)\geq 1$ which implies that $HF_\mathrm{red}(M,\mathfrak{s}_0)$ contains an odd number of cyclic summands whereas $HF_\mathrm{red}(M,\mathfrak{s}_1)$ contains an even number of cyclic summands. So $HF_\mathrm{red}(M,\mathfrak{s}_0)\not\cong HF_\mathrm{red}(M,\mathfrak{s}_1)$.

Suppose that $\mathfrak{s}_i$ is defined for $i=0$ or $1$. It follows from the action of $\tau$ that the set $S_i$ of spin$^c$-structures
\[
S_i=\{ \mathfrak{s}\in \mathrm{Spin}^c(M)\mid HF_\mathrm{red}(M,\mathfrak{s})\cong HF_\mathrm{red}(M,\mathfrak{s}_i)\}
\]
has odd cardinality. Since Heegaard Floer is invariant under conjugation, it follows that $S_i$ contains a self-conjugate spin$^c$-structure. Since $M$ has one self-conjugate  spin$^c$-structure if $p$ is odd and two if $p$ is even, it follows that the reduced Heegaard Floer groups in the self-conjugate spin$^c$-structures are isomorphic to $HF_\mathrm{red}(M,\mathfrak{s}_0)$ and $HF_\mathrm{red}(M,\mathfrak{s}_1)$, according to their existence.  
\end{proof}
For $K$ a non-trivial $L$-space knot of genus $g(K)$, the non-zero coefficients of $\Delta_K(t)$ are equal to $\pm 1$ and alternate in sign. That is the symmetrized Alexander polynomial can be written in the form
\[
\Delta_K(t)=\sum_{k=-\ell}^\ell (-1)^{k+\ell}t^{\beta_k}
\]
for some sequence of integers $\beta_l >\dots > \beta_{-l}$ such that $\beta_k=-\beta_{-k}$, $\beta_\ell=g(K)$ and $\beta_{\ell-1}=g(K)-1$. From such a description, one sees that the torsion coefficients form a non-increasing sequence such that $t_g(K)=0$ and $t_{g-1}=1$. It follows that if one takes $T_1$ to be the quantity
\[T_1=|\{i\mid t_i(K)=1 \}|,\]
then $t_i(K)=1$ if and only if
\[g-T_1\leq i \leq g-1.\]

\begin{theorem}
\label{thm:neg_L_space_obstruction}
    Let $M$ be obtained by $-p/q$-surgery on an $L$-space knot $K$ for $p/q>0$. Then $M$ admits a co-orientable taut foliation. However, if $M$ admits a co-orientable taut foliation with Euler class zero, then
    \[
    q(2g(K)-1) = r + mp 
    \]
    for some integers $m \geq 0$ and $1\leq r< 2q T_1$, where $T_1=|\{i\mid t_i(K)=1 \}|$.
\end{theorem}

\begin{proof}
    An $L$-space knot is fibred and its monodromy has positive fractional Dehn twist coefficient \cite{Ni07, Hed10}. Consequently Roberts \cite{RobertsSurfacebundle2} shows that any negative surgery on an $L$-space knot admits a co-orientable taut foliation.

    Now suppose that $M$ admits a co-orientable taut foliation with Euler class zero.
    From Lemma~\ref{lem:negative_surgery_calc}, Proposition~\ref{prop:euler_class_obstruction} and the discussion prior to this proposition, we see that there is a non-negative integer $m$ such that $pm+q-1$ is even and
    \[
    g-T_1\leq \left\lfloor \frac{pm+q-1}{2q} \right\rfloor \leq g-1.
    \]
    Since $pm+q-1$ is even, $\frac{pm+q-1}{2q}$ can be expressed as a fraction with denominator $q$ and so we have
    \[
    g-T_1\leq \frac{pm+q-1}{2q} \leq g-\frac{1}{q}.
    \]
    Rearranging this implies that
    \[
    q(2g-1)-2qT_1+1\leq pm \leq q(2g-1)-1.
    \]
    This implies that there is $1\leq r\leq 2q T_1 -1$ such that
    \[
    q(2g-1) = r + mp.
    \]
\end{proof}

\subsection{Applications of Theorem~\ref{thm:neg_L_space_obstruction}}
Using Theorem~\ref{thm:neg_L_space_obstruction}, we can generate many examples of manifolds admitting a co-orientable taut foliation, but none with Euler class zero. For instance, it allows us to drop the condition that $q$ be even from Theorem~\ref{thm: p/q-surgery} for many $L$-space knots.

\begin{corollary}
\label{cor: q arbitrary t0}
    Let $M$ be a manifold obtained by $-p/q$-surgery on an $L$-space knot $K$ with $t_0(K) \geq 2$, where $p/q \geq 2g(K) -1$. Then $M$ admits a co-orientable taut foliation but none with Euler class zero.  
\end{corollary}

\begin{proof}
If $M$ admits a co-orientable taut foliation with zero Euler class then Theorem~\ref{thm:neg_L_space_obstruction} implies that 
$q(2g(K)-1) = r + mp$ for integers $m \geq 0$ and $r$, where $1\leq r < 2q T_1$. If $m > 0$, then $p/q = (2g(K)-1)/m - r/qm < 2g(K) -1$, contrary to our hypotheses. Thus $m = 0$ and therefore $q(2g(K)-1) = r < 2qT_1$. On the other hand, the condition that $t_0(K) \geq 2$ implies that $T_1 \leq g(K) - 1$, so $q(2g(K)-1) = r < 2q(g(K) - 1)$, which is impossible. Thus $M$ admits no co-orientable taut foliation with Euler class zero.  
\end{proof}

\begin{remark}
\label{rem: t_0 condition}
The condition $t_0(K) \geq 2$ appears to be generic among $L$-space knots. The authors are aware of only five $L$-space knots for which it does not hold: the three torus knots $T(2, 3)$, $T(2, 5)$ and $T(3, 4)$, the $(2, 3)$-cable of $T(2, 3)$, and the $(2,5)$-cable of $T(2, 3)$, of genus $1,2,3,3,$ and $4$, respectively. See \cite{Krcatovich18}. Since $t_0(K) = 1$ if and only if $\Delta_K(t) = (t^g + t^{-g}) - (t^{g-1} + t^{-(g-1)}) + 1$, where $g = g(K)$, it follows from \cite[Corollary 4.5]{Tange20} that if $K$ has a lens space surgery (so, conjecturally, is a Berge knot) and $g(K) \ge 4$ then $t_0(K) \ge 2$. 
\end{remark}

It will be convenient to note that one can bound the quantity $T_1$ appearing in Theorem~\ref{thm:neg_L_space_obstruction} using the bridge index.
\begin{corollary}\label{cor:bridge_bound}
Let $K$ be an $L$-space knot of bridge index $\mathrm{br}(K)$ and let $M$ be obtained by $-p/q$-surgery on $K$ for $p/q>0$. Then $M$ admits a co-orientable taut foliation.  However, if $M$ admits a co-orientable taut foliation with Euler class zero, then
\[
q(2g(K)-1) = r + mp 
\]
for some integers $m \geq 0$ and $1\leq r< 2q \mathrm{br}(K)$.
\end{corollary}

\begin{proof}
Given Theorem~\ref{thm:neg_L_space_obstruction}, it suffices to prove that $T_1\leq \mathrm{br}(K)$. This follows from the properties of the invariant $\mathrm{Ord}_v(K)$ defined in \cite{JMZ2020}. The symmetrized Alexander polynomial of a non-trivial $L$-space knot can be written in the form
    \[
    \Delta_K(t)=t^g-t^{g-1}+ t^{g-T_1}+ \text{lower order terms}
    \]
    The calculation of $\mathrm{Ord}_v(K)$ for $L$-space knots \cite[Lemma~5.1]{JMZ2020} shows that
    \[T_1\leq \mathrm{Ord}_v(K)+1.\]
    However, for any knot $K$ the bound $\mathrm{Ord}_v(K)+1\leq \mathrm{br}(K)$ holds \cite[Corollary~1.9]{JMZ2020}.
\end{proof}

\begin{customthm}{\ref{thm:hyperbolic_examples}}
    For any integer $p\geq 7$ there are infinitely many hyperbolic 3-manifolds $M$ such that $H^2(M)\cong \mathbb{Z}/p$ and $M$ admits a co-orientable taut foliation, but none with Euler class zero. Further, $\pi_1(M)$ is left-orderable.
\end{customthm}
 
\begin{proof}
    For $g\geq 5$, let $K_g$ denote the $(-2,3,2g-3)$-pretzel knot. The knot $K_g$ is a hyperbolic $L$-space knot \cite{Ozsvath2005lens}. Moreover, we have $br(K_g)=3$ and $g(K)=g$. Let $M$ be obtained by $-p/q$-surgery on $K_g$ for $p/q>0$.  The surgery slope $4g-1$ is an exceptional surgery slope for $K_g$ since it yields a Seifert fibred space \cite[Proposition~16]{BleilerHodgson}.  Since $\Delta(4g-1, -p/q)=|(4g-1)q+p|\geq 20$, the known bounds on exceptional fillings imply that $M$ is a hyperbolic manifold \cite{LaMe13}. The work of Roberts shows that there exists a co-orientable taut foliation on $M$. Moreover, if $g\neq 6$, then $M$ has left-orderable fundamental group \cite{Tran, Var21}.
    
    On the other hand, if $M$ admits a co-orientable taut foliation with Euler class zero, then Corollary~\ref{cor:bridge_bound} implies
    \[
    q(2g(K_g)-1)\equiv r \; (\bmod \; p)
    \]
    where $1\leq r\leq 6q-1$. Thus we see that $M$ does not admit a co-orientable taut foliation with Euler class zero, whenever $p/q\geq 7$ and $g\equiv 0 \; (\bmod \; p)$. 
\end{proof}

\section{The Euler class of the normal bundle of a Seifert fibration} 
\label{sec: Euler class normal bundle}

Let $M$ be a closed, oriented Seifert fibred manifold whose base orbifold has an 
orientable underlying space. In this section, we compute the Euler class of  $\nu_M$, the normal bundle to the Seifert fibres of $M$, and determine conditions for it to vanish. 

\subsection{Seifert notation}
\label{subsection: notation}

We first recall the Seifert notation that we will use to describe $M$. 

Consider $M$ constructed by attaching $n+1$ solid tori $N_i$, $i=0, 1, \dots,n$, 
to a trivial $S^1$-bundle over an orientable surface $F$ with $|\partial F| = n+1$, 
where we assume that $N_0$ is a tubular neighbourhood of a regular fibre and the 
$N_i$ are tubular neighbourhoods of the singular fibres for $i=1, 2, \ldots,n$. 
Denote this trivial $S^1$-bundle by $P \cong F \times S^1$. 
Fix orientations on $F$ and on the fibres so that they intersect  in $M$.

Let $$\partial F = \mu_0 \sqcup \mu_1 \sqcup \cdots \sqcup \mu_n$$ with the induced orientation, and  
\[
\partial P = T_0 \sqcup T_1 \sqcup \cdots \sqcup T_n
\]
with $\mu_i \subset T_i$ for each $i = 0, 1, \cdots, n$.  So on each $T_i$, we use the basis given by $\{h_i, \mu_i\}$, where $h_i$ is given by a Seifert fibre.  By our orientation convention, we have $\mu_i\cdot h_i = 1$. 

Let $D_i$ be a meridional disk of $N_i$, oriented so that it intersects the Seifert fibres 
positively. Under the attaching map, we have $\partial D_i = - \bigl(a_i h_i + b_i \mu_i\bigr)$ 
as oriented curves, for some  $a_i \in \mathbb{Z}$, $b_i \in \mathbb{Z}_{>0}$ and $a_i$, $b_i$ coprime. In fact, since $N_0$ is a tubular neighbourhood of a regular fibre, and the $N_i$ are tubular neighbourhoods of singular fibres, we have $b_0 = 1$ and $b_i\geq 2$. 

Let $g\geq 0$ denote the genus of $F$, $b = \frac{a_0}{b_0} \in \mathbb{Z}$, and denote 
\[
M = M\Bigl(g;\, b, \tfrac{a_1}{b_1}, \ldots, \tfrac{a_n}{b_n}\Bigr)
\]

Let $\mathcal{B}$ denote the base orbifold of $M$, 
which is a genus $g$ closed orientable surface with $n$ cone points of orders $b_1,\ldots,b_n$ respectively. We use $|\mathcal{B}|$ to denote the underlying space of $\mathcal{B}$.

\begin{remarks}
    \label{remark: Seifert notation}
    
\begin{enumerate}[leftmargin=*]
    \setlength{\itemsep}{1em}
\item This description of the Seifert fibred manifold is not unique. More precisely,  
$M\Bigl(g;\, \tfrac{a_0}{b_0}, \ldots, \tfrac{a_n}{b_n}\Bigr)$ is homeomorphic to $M\Bigl(g;\, \tfrac{a'_0}{b_0}, \ldots, \tfrac{a'_n}{b_n}\Bigr)$ by an orientation-preserving homeomorphism preserving fibres if and only if up to reordering, 
\[
\frac{a_i}{b_i} - \frac{a'_i}{b_i} \in \mathbb{Z}
\quad \text{and} \quad 
\sum_{i=0}^n \frac{a_i}{b_i} = \sum_{i=0}^n \frac{a'_i}{b_i}.
\]
By choosing exactly one of the $N_i$, namely $N_0$, to be a neighbourhood of a regular fibre in the construction above, we are adopting the usual convention of having only one of the $\tfrac{a_i}{b_i}$ be an integer.

\item 
Given $M = M\bigl(g;\, b, \tfrac{a_1}{b_1}, \ldots, \tfrac{a_n}{b_n}\bigr)$, the sum $$e(M) := b + \sum_{i=1}^n \frac{a_i}{b_i}$$
is an invariant of the given Seifert fibring of $M$ called the \emph{Euler number}. 

\item The effect on $b, a_1, \ldots , a_n$ and $e(M)$ of reversing the orientation of $M$ is to replace them by their negatives.

\end{enumerate}
\end{remarks}

\begin{remark}
    \label{rem: n =0 }When $n=0$, the manifold $M = M(g; b)$ is an $S^1$-bundle over a closed, 
orientable surface $B$ of genus $g$ with Euler number $b$. Moreover, the necessary and 
sufficient condition in Theorem~\ref{thm: Euler class zero} simplifies to the existence of 
$m \in \mathbb{Z}$ such that 
\begin{equation}
    \label{eqn: S1 bundle condition}
mb = \chi(B)
\end{equation}
By assumption, the Euler class of the circle bundle $f: M \to B$ is $bu$, where $u$ is a generator of $H^2(B) \cong \mathbb Z$. It is easy to see that $\nu_M$ is the pull-back by $f$ of the tangent bundle $\tau_B$ of $B$. Hence, 
\begin{equation}
    \label{eqn: enu} 
    e(\nu_M) = f^*(e(\tau_B)), 
    \end{equation} 
Let $\langle h \rangle$ denote the central subgroup of $\pi_1(M)$ generated by a circle fibre of $f$ and consider the $5$-term exact sequence in cohomology associated to the central extension $1 \to \langle h \rangle \to \pi_1(M) \to \pi_1(B) \to 1$:
\[0 \to H^1(B)\xrightarrow{} H^1(M) \to \mathbb Z \xrightarrow{\; \delta \;} H^2(B) \xrightarrow{\; f^* \;} H^2(M) 
\]
It follows from \cite[Theorem 4]{HS53} that the image of $-1$ under the transgression homomorphism $\delta$ is the Euler class $bu$ of the circle bundle. Hence $\mbox{kernel}(f^*)$ is generated by $bu$. Then by (\ref{eqn: enu}), $e(\nu_M) = 0$ if and only $e(\tau_B) = \pm \chi(B)u$ is a multiple of $bu$, which is condition  (\ref{eqn: S1 bundle condition}). 

We consider the following two special cases: 
\begin{enumerate}
    \item If $M = M(g; b)$ is a rational homology sphere, then the base $B$ must be the $2$-sphere, 
    and $M = M(0; b)$ is a lens space $L(|b|, 1)$. 
    In this case, the Euler class $e(\nu_M)$ vanishes if and only if $b \in \{\pm 1, \pm 2\}$. 
    
    \vspace{.2cm} \item When $B$ is the torus, (\ref{eqn: enu}) implies that the Euler class $e(\nu_M)$ is always trivial, since the tangent bundle $T\mathcal{B}$ is trivial, 
   
\end{enumerate}

\end{remark}

\subsection{The triviality of the normal bundle of a Seifert fibring} 
Here is our main result.

\begin{theorem}
\label{thm: Euler class zero}
Given a closed oriented Seifert fibred manifold $M = M\bigl(g;\, b, \tfrac{a_1}{b_1}, \ldots, \tfrac{a_n}{b_n}\bigr)$, the Euler class of the 
normal bundle $\nu_M$ of the Seifert fibration vanishes if and only if there exists $m \in \mathbb{Z}$ such that
\vspace{-.2cm}
\begin{enumerate}
\vspace{.2cm} \item[\rm (1)] $m a_i \equiv  1 \, (\text{{\rm mod} }b_i), \text{ for } i =  1, \cdots, n, \; \text{and}$ 
\vspace{.2cm} \item[\rm (2)] $ me(M) = \chi(\mathcal{B})$.

\end{enumerate}

\end{theorem}

\begin{corollary}
    \label{cor: necessary condition}
Suppose that $M$ is a closed oriented Seifert fibred manifold $M$ such that the Euler class of $\nu_M$ vanishes. Then one of the following holds:
\vspace{-.2cm}
\begin{enumerate}
\vspace{.2cm} \item[\rm (1)]   $M$ is an orientable $S^1$-bundle over the torus, 
\vspace{.2cm} \item[\rm (2)]  $e(M) = \chi(\mathcal{B}) = 0$, 
\vspace{.2cm} \item[\rm (3)] $ e(M)\neq 0, \, \chi(\mathcal{B}) \neq 0, $ and $$\frac{\chi(\mathcal{B})}{e(M)}\in \mathbb{Z}.$$ 

\end{enumerate}
 
\end{corollary}

The following lemma is key to applying these results to the Euler classes of foliations on Seifert fibred rational homology $3$-spheres. 

\begin{lemma}
\label{lem: same Euler class}
Suppose that $M$ is a Seifert fibred rational homology $3$-sphere. Then 

$(1)$ Any co-orientable taut foliation on $M$ is horizontal. That is, it can be isotoped to be everywhere transverse to the Seifert fibres of $M$.

$(2)$ The tangent space to any co-orientable taut foliation on $M$ is isomorphic to $\nu_M$ as an unoriented bundle. In particular, the Euler class of any co-oriented taut foliation on $M$ agrees, up to sign, with that of $\nu_M$.   
\end{lemma}

\begin{proof}
That any co-oriented taut foliation on $M$ must be horizontal is \cite[Corollary 7]{Brittenham93} in the case $n = 3$ and \cite[Proposition 6.1]{BC17} in general. Thus (1) holds. Consequently, the tangent bundle of any such foliation is isomorphic to $\nu_M$ (as an unoriented bundle), which implies (2). 
\end{proof}

Theorem \ref{thm: main application} follows immediately from this lemma and Theorem \ref{thm: Euler class zero}. 

\subsection{Proof of Corollary \ref{cor: necessary condition}}
\label{subsec: condns}
We prove Corollary~\ref{cor: necessary condition} using Theorem \ref{thm: Euler class zero} here and defer the proof of 
Theorem~\ref{thm: Euler class zero} to the next subsection. 

Let $\nu_M$ denote the normal plane bundle over $M$ to the Seifert fibration, oriented with the induced orientation from those on $M$ and on the fibres.

\begin{proof}[Proof of Corollary \ref{cor: necessary condition}]
Assume that the Euler class of the normal bundle vanishes. 

We first consider the case $n=0$, in which $\mathcal{B}$ is a genus $g$ surface and 
$M$ is an $S^1$-bundle over $\mathcal{B}$ of Euler number $b$. The necessary and sufficient condition in 
Theorem \ref{thm: Euler class zero} is equivalent to that the existence of $m \in \mathbb{Z}$ 
such that $mb = \chi(\mathcal{B})$. 

If $m = 0$, then $\chi(\mathcal{B}) = mb = 0$ and therefore $M$ is an oriented
$S^1$-bundle over the torus. On the other hand, if $M$ is not an oriented
$S^1$-bundle over the torus, then $\chi(\mathcal{B}) \neq 0$. Then for the equality 
$mb = \chi(\mathcal{B})$ to hold, we must have $e(M) = b \neq 0$. Hence 
\[
m = \frac{\chi(\mathcal{B})}{e(M)} \in \mathbb{Z}.
\]

Next we consider the case $n > 0$. Again by Theorem~\ref{thm: Euler class zero}, there exists $m \in \mathbb{Z}$ such that 
$m a_i \equiv 1 \pmod{b_i}$ and  $m e(M) = \chi(\mathcal{B})$. Since $m$ cannot be zero in this case, either (1) or (2) holds. 
\end{proof}

\subsection{Proof of Theorem \ref{thm: Euler class zero}}
\label{subsec: proof main theorem}
Here we provide a proof of Theorem~\ref{thm: Euler class zero}.
When the base orbifold is hyperbolic, i.e., $\chi(\mathcal{B}) < 0$, we will give another 
perspective on Theorem~\ref{thm: Euler class zero} and Corollary \ref{cor: necessary condition} in \S \ref{sec: normal}. There we will show that the Euler class of the normal bundle is precisely the pullback of the Euler class of a discrete faithful representation 
of $\pi_1(\mathcal{B})$ in $PSL(2,\mathbb{R})$ under the homomorphism $\pi_1(M) \to \pi_1(\mathcal{B})$ induced by the natural projection $M\rightarrow \mathcal{B}$ \footnote{Since $M$ is irreducible when $\mathcal{B}$ is hyperbolic, 
we have $H^2(M) \equiv H^2(\pi_1(M))$.}.

\begin{proof}[Proof of Theorem \ref{thm: Euler class zero}]
To simplify the notation in the proof, we will use $\tfrac{a_0}{b_0}$ in place of $b$, 
as explained in \S\ref{subsection: notation}. Thus
\[
M = M\bigl(0;\, b, \tfrac{a_1}{b_1}, \ldots, \tfrac{a_n}{b_n}\bigr) 
   = M\bigl(0;\, \tfrac{a_0}{b_0}, \tfrac{a_1}{b_1}, \ldots, \tfrac{a_n}{b_n}\bigr)
\]
We will compute the Euler class of the normal bundle, denoted by $e(\nu_M) \in H^2(M)$, using the following exact sequence of the pair $(M, \sqcup_i T_i)$:
\[
\cdots \to \oplus_i H^1(T_i) \xrightarrow{\delta} H^2(M, \sqcup_i T_i) \xrightarrow{i^*} H^2(M) \to 0
\]

Note that 
\[
H^2(M, \sqcup_i T_i) \cong H^2(P, \sqcup_i T_i) \oplus H^2(N_1, T_1) \oplus \cdots \oplus H^2(N_n, T_n),
\]
and by the K\"unneth formula,
\[
H^2(P, \sqcup_i T_i) \cong H^1(F, \partial F) \otimes H^1(S^1) \oplus H^2(F,\partial F)
\]
Recall that $\partial F = \sqcup_i \mu_i$ with $\mu_i \subset T_i$ for each $i = 0,1, \ldots, n$.  Suppose that oriented simple closed curves $\alpha_r,\beta_r$, $r= 1, \ldots, g$, together with oriented arcs $c_j$ on $F$ that connect the boundary component $\mu_j$ to the boundary component $\mu_{j+1}$, $j=0, 1,\ldots,n-1$ when $n>0$, represent a set of generators of $H_1(F,\partial F)$. Then $H_2(P,\sqcup_i T_i)$ is freely generated by the classes of vertical tori $\alpha_r \times S^1$, $\beta_r \times S^1$, the horizontal surface $F$ and vertical annuli $c_j \times S^1$ when $n>0$. We denote their duals in $H^2(P, \partial P)$ by 
\[
\alpha_r^* \otimes S^*, \quad \beta_r^* \otimes S^*, \quad c_j^* \otimes S^*, \quad F^\ast
\]
respectively, where  $r=1,\ldots,g$ and $j=0, 1,\ldots,n-1$.

Finally, let  $D_i^* \in H^2(N_i, T_i)$ denote the dual of $[D_i,\partial D_i] \in H_2(N_i,T_i)$ for each $i = 0, 1, \ldots, n$.  

In summary, a cohomology class in
\[
H^2(M, \sqcup T_i) \cong H^1(F,\partial F) \otimes H^1(S^1) \oplus H^2(F,\partial F) \oplus H^2(N_1, T_1) \oplus \cdots \oplus H^2(N_n, T_n)
\]
can be uniquely expressed as a linear combination of
\[
\alpha_r^* \otimes S^*, \quad \beta_r^* \otimes S^*, \quad F^*, \quad D_i^* \quad \text{and} \quad c_j^* \otimes S^* \text{ when } n> 0
\]
where $r=1,\ldots,g$, $j=0, 1,\ldots,n-1$ and $i=0, 1,\ldots,n$. 

The proof proceeds in two steps (1) and (2). 
\begin{enumerate}[leftmargin=*]
    \setlength{\itemsep}{1em}
\item Find the preimage of the Euler class $e(\nu_M)$ of the normal bundle $\nu_M$ in $H^2(M, \sqcup T_i)$.

\medskip 

\noindent
Let $\sigma$ be a nowhere vanishing vector field along $\partial P = \sqcup T_i$ that points out of $P$. Since $T_i$ is tangent to the Seifert fibres, $\sigma$ is a section of  $\nu_M$ along $\sqcup T_i$. We give a cellular structure on $M$ for which we assume that 
the vertical surfaces $\alpha_r \times S^1$, $\beta_r \times S^1$, $c_j \times S^1$, and the horizontal surfaces $F$ and $D_i$, as well as the $T_i$, all 
belong to the 2-skeleton. We also assume that the curves $\alpha_r$, $\beta_r$ and the arcs $c_j$ on $F$ are in the 1-skeleton.  

From $\sigma$ one can define an obstruction $2$-cochain, denoted by $c_\sigma$, that represents 
the relative Euler class of the normal bundle $\nu_M$ in $H^2(M,\sqcup T_i)$. By construction, 
the image of $[c_\sigma]$ under $H^2(M,\sqcup T_i) \to H^2(M)$ is the Euler class $e(\nu_M)$. (See \cite[Section 2]{Hu23} for instance.)

To write $[c_\sigma]$ as a linear combination of 
$\alpha_r^* \otimes S^*$, $\beta_r^* \otimes S^*$, $c_j^* \otimes S^*$, $F^*$, and $D_i^*$, 
it suffices to compute its evaluation on the dual homology classes.

\medskip

\begin{enumerate}
    \setlength{\itemsep}{1em}
\item $[c_\sigma]$ evaluates to zero on all vertical surfaces.  
This is because, starting from $\sigma|_{\partial F}$, where $\partial F  \subset \partial P = \sqcup T_i$, 
there is no obstruction to extending $\sigma$ to the 1-skeleton of $F$. Since 
$P \cong F \times S^1$ is a trivial $S^1$-bundle, one can extend the section further to 
the vertical surfaces $\alpha_r\times S^1$, $\beta_r\times S^1$ and $c_j\times S^1$ by staying constant along the fibre direction.  

\item By the Poincar\'e-Hopf Theorem, $[c_\sigma]$ evaluated on a horizontal surface equals its Euler characteristic, i.e.,
\[
[c_\sigma]([D_i,\partial D_i]) = 1, \quad i=0, 1,\ldots,n, \]
and \[
[c_\sigma]([F,\partial F]) = \chi(F) = \chi(|\mathcal{B}|) - (n+1).
\]
\end{enumerate}

\item Compute the image of $\delta: \oplus_i H^1(T_i) \rightarrow H^2(M, \sqcup_i T_i)$. 

\medskip

\noindent
Let $h_i^*, \mu_i^* \in H^1(T_i)$ denote the ($\mbox{Hom}$) dual classes of $[h_i], [\mu_i] \in H_1(T_i)$ 
for $i=0, 1,\ldots,n$. We compute $\delta h_i^*$ and $\delta \mu_i^*$, which generate the 
image of $\delta$.  

\medskip 
\begin{enumerate}
    \setlength{\itemsep}{1em}
\item It is easy to see that, since $\delta h_i^* = h_i^* \circ \partial$ and 
$\delta \mu_i^* = \mu_i^* \circ \partial$, they vanish on all vertical tori, 
$\alpha_r \times S^1$ and $\beta_r \times S^1$.

 \item As $\partial([c_j \times S^1]) = [h_{j+1}] - [h_j]$ in $H_1(\sqcup_i T_i)$ for $j=0, 1,\ldots,n-1$,  
we have
\[
\delta \mu_i^*([c_j \times S^1]) = 0 \quad \text{for all } i=0, 1,\ldots,n, \text{ and } \ j= 0, 1,\ldots,n-1,
\]
and
\[
\delta h_i^*([c_j \times S^1]) =
\begin{cases}
1 & \text{if } j= i - 1 \\[6pt]
-1 & \text{if } j = i \\[6pt]
0 & \text{otherwise}
\end{cases}
\]

\item Since $\partial F = \sqcup_i \mu_i$, we have $\delta h_i^*([F,\partial F]) = 0$ and 
$\delta \mu_i^*([F,\partial F]) = 1$ for all $i=0, 1,\ldots,n.$

\item Finally, to compute the values of  $\delta h_i^*$ and $\delta \mu_i^*$ over $[D_i,\partial D_i] \in H_2(N_i,T_i)$, we note that  $\partial [D_i] = -(a_i [h_i] + b_i [\mu_i])$ by our orientation convention in order to keep $b_i > 0$ (see \S \ref{subsection: notation}). Then it follows that for $i = 0, 1, \ldots, n$
\[
\delta h_i^*([D_i,\partial D_i]) = h_i^*(-a_i [h_i] - b_i [\mu_i]) = -a_i
\]
and
\[
\delta \mu_i^*([D_i,\partial D_i]) = \mu_i^*(-a_i [h_i] - b_i [\mu_i]) = - b_i.
\]
\end{enumerate}
\end{enumerate}
In summary, we have 
\[
[c_\sigma] = \chi(F) F^* + \sum_{i=0}^n D_i^*
\]
The image of $\delta$ is generated by
\[
\delta \mu_i^* = F^* - \sum_{i=0}^n b_i D_i^*, \quad \text{ and } \quad
\delta h_i^* =
\begin{cases}
 -c_0^* \otimes S^* - a_0 D_0^*, & \text{for } i=0, \\[6pt]
  c_{i-1}^* \otimes S^* - c_i^* \otimes S^* - a_i D_i^*, & \text{for } 0<i<n, \\[6pt]
  c_{n-1}^* \otimes S^* - a_n D_n^*, & \text{for } i=n.
\end{cases}
\]

if $n>0$; and when $n=0$, the image of $\delta$ is generated by
\[
\delta \mu_0^* = F^* - b_0 D_0^* = F^* - D_0^* 
\quad \text{and} \quad 
\delta h_0^* = -a_0 D_0^* = -b D_0^*.
\]
The Euler class $e(\nu_M)=i^*([c_\sigma]) \in H^2(M)$ vanishes 
if and only if $[c_\sigma] \in \mathrm{Im}(\delta)$. That is, there exist integers 
$m_i, n_i$, for $i=0, 1,\ldots,n$, such that
\[
\delta\Bigl(\sum_i -m_i h_i^* + \sum_i n_i \mu_i^* \Bigr) = [c_\sigma].
\]

\begin{enumerate}[leftmargin=*]
    \setlength{\itemsep}{0.5em}
    \item When $n>0$, for the coefficients of $c_j^* \otimes S^*$ to be zero, the $m_i$ must be identical. 
Hence it follows that $[c_\sigma] \in \mathrm{Im}(\delta)$ if and only if there exist 
$m, n_i \in \mathbb{Z}$ for $i = 0, \ldots, n$ such that
\[
\begin{cases}
 m a_i - n_i b_i = 1, \mbox{ for } i = 0, 1, \ldots, n \\[6pt]
 \sum_{i=0}^n n_i = \chi(F) = \chi(|\mathcal{B}|) - (n+1).
\end{cases}
\]

Recall that $\tfrac{a_0}{b_0} = b$. Hence the equation 
$m a_0 - n_0 b_0 = 1$ simplifies to $m b - n_0 = 1$, and therefore, the condition above 
is equivalent to the existence of $m, n_i \in \mathbb{Z}$ for $i=1,\ldots,n$ such that 
\[
\begin{cases}
 m a_i - n_i b_i = 1, \text{ for } i = 1,\ldots,n \\[6pt]
 mb - 1 + \displaystyle\sum_{i=1}^n n_i = \chi(F) = \chi(|\mathcal{B}|) - (n+1)
\end{cases}
\]

For $i \geq 1$, we have 
\[
n_i = \frac{m a_i - 1}{b_i} 
    = \frac{m a_i}{b_i} - \frac{1}{b_i} 
\]
Hence, the second condition is equivalent to 
\[m(b + \sum_{i=1}^n \frac{a_i}{b_i}) = \chi(|\mathcal{B}|) - n +  \sum_{i=1}^n \frac{1}{b_i}\]
Note that 
\[
e(M) = b + \sum_{i=1}^n \frac{a_i}{b_i}
\quad \text{and} \quad
\chi(\mathcal{B}) = \chi(|\mathcal{B}|) - n + \sum_{i=1}^n \frac{1}{b_i}
\]
This completes the proof of the theorem when $n>0$. 

\item 
If $n=0$ (and hence there is no singular fibre), the condition above reduces to the existence of  $m\in \mathbb{Z}$ such that 
\[mb = \chi(\mathcal{B})\] 
From the discussion above, one can easily verify that this is a necessary and 
sufficient condition for the Euler class to vanish in this special case.

\end{enumerate}

This proves the theorem.
\end{proof}

\section{Representations of Fuchsian groups and the  normal bundle of a Seifert fibration}
\label{sec: normal} 

In this section we provide alternative proofs of Theorem \ref{thm: Euler class zero} and Corollary \ref{cor: necessary condition} from the perspective of discrete faithful representations of Fuchsian groups. We normalise our notation for oriented Seifert manifolds 
$M = M\bigl(g;\, b, \tfrac{a_1}{b_1}, \ldots, \tfrac{a_n}{b_n}\bigr)$ 
by further requiring that $\tfrac{a_i}{b_i} \in (0,1)$ for each $i=1,\ldots,n$. See Remark \ref{remark: Seifert notation}.

\begin{figure}[ht]
\includegraphics[scale=1.4]{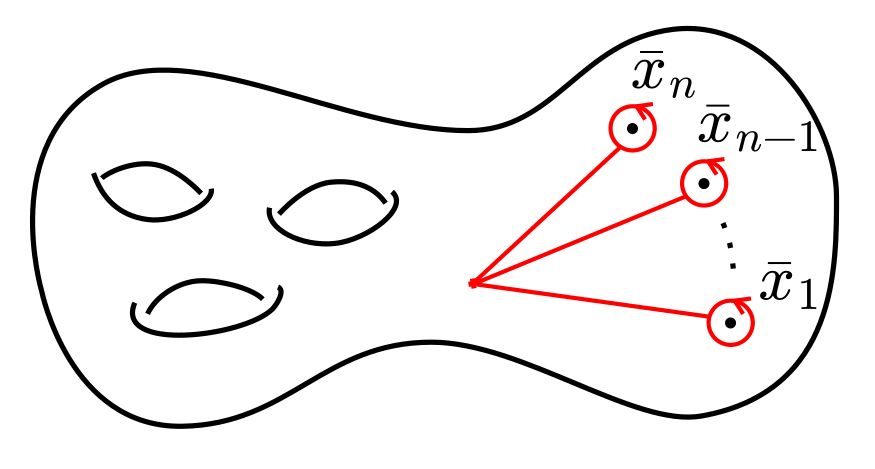}
\caption{The generators $\bar{x}_i$ in the presentation of $\Gamma = \pi_1(\mathcal{B})$}
\label{fig: x1}
\end{figure}
The fundamental group of $M$ has a presentation
\begin{eqnarray} 
\label{eqn: pres 2} 
\hspace{1cm} \langle x_1, x_2, \ldots, x_n, w_1, w_2, \ldots, w_g,  z_1, z_2, \ldots , z_g, h \;|\; h \mbox{ central}, x_i ^{b_i} = h^{-a_i}  \mbox{ for each } i, \\ 
\hspace{2cm} [w_1, z_1][w_2, z_2] \cdots [w_g, z_g] x_1x_2 \cdots x_n = h^{b} \rangle \nonumber, 
 \end{eqnarray} 
and $\Gamma = \pi_1(\mathcal{B})$ has a presentation 
\begin{eqnarray} 
\label{eqn: pres 1} 
\hspace{1cm} \langle \bar x_1, \ldots, \bar x_n, \bar w_1, \ldots, \bar w_g,  \bar z_1,  \ldots , \bar z_g \;|\; \bar x_i ^{b_i} = 1,  [\bar{w}_1, \bar{z}_1] \cdots [\bar{w}_g, \bar{z}_g]\bar x_1 \cdots \bar x_n = 1 \rangle 
\end{eqnarray}
These two presentations are related via a well-known central extension: 
\begin{equation} 
\label{eqn:  central extension} 
E_M: 1 \to \mathbb Z \to \pi_1(M) \xrightarrow{\; \alpha \;} \Gamma \to 1, 
\end{equation}
Here $h$ corresponds to the image in $\pi_1(M)$ of $1 \in \mathbb Z$, $\alpha(x_i) = \bar x_i$ for each $i$, $\alpha(w_j) = \bar w_j$ for each $j$ and $\alpha(z_k) = \bar z_k$ for each $k$. 

\subsection{The second cohomology group of $\Gamma$} 
\label{subsec: H^2}
Jankins and Neuman determined a canonical isomorphism between $H^2(\Gamma)$ and the abelian group having generators $A_0, A_1, \ldots, A_n$ and relations $A_0 = b_i A_i$ for $1 \leq i \leq n$ \cite{JN85a}. (The case that $g = 0$ is discussed in \cite{JN85a} and the general case in \cite{JN83}.) Each element of $H^2(\Gamma)$ can be uniquely expressed in the form
\begin{equation} 
\label{eqn: normal form} 
dA_0 + c_1A_1 + \cdots + c_n A_n, \nonumber
\end{equation}
where  $d, c_1, \ldots, c_n \in \mathbb Z$ and $1 \leq c_i \leq b_i - 1$ for each value of $i$.

The Euler class $e(E_M) \in H^2(\Gamma)$ is the cohomology class which classifies the extension (\ref{eqn:  central extension}). In \cite{JN85a} it was shown to be 
\begin{equation} 
\label{eqn: e(E_M)}
e(E_M) = -(bA_0 + a_1A_1 + a_2 A_2 + \cdots + a_n A_n)
\end{equation}
This class arises as the image of $-1$ under the transgression homomorphism $\delta$ from the $5$-term exact sequence in cohomology of $E_M$:
\begin{equation}
\label{eqn: 5-term}
0 \to H^1(\Gamma)\xrightarrow{} H^1(M) \to \mathbb Z \xrightarrow{\; \delta \;} H^2(\Gamma) \xrightarrow{\; \alpha^* \;} H^2(M) 
\end{equation}
See \cite[Theorem 4]{HS53}. Thus, 
\begin{equation}
\label{eqn: kernel alpha*}
\mbox{kernel}(\alpha^*) = \langle e(E_M) \rangle
\end{equation}

 \subsection{The normal bundle of the Seifert fibring of $M$ and discrete representations} 
Since $\mathcal{B}$ is hyperbolic, its universal cover can be identified with $\mathbb H^2$ and the associated group of deck transformations with a discrete subgroup of $\mbox{Isom}_+(\mathbb H^2) = PSL_2(\mathbb R)$. This yields a discrete faithful representation $\bar \varphi: \Gamma \to PSL_2(\mathbb R)$. 

There is a Seifert fibring $\widehat M \to \mathbb H^2$ lifting $M \to \mathcal{B}$ for which the following diagram commutes
\begin{center} 
\begin{tikzpicture}[scale=0.8]
 
\node at (9, 4.5) {$\widehat M = \mathbb H^2 \times S^1$};
\node at (15.1, 4.5) {$\mathbb H^2$};
 
\node at (9, 2) {$M$};
\node at (15.1, 2) {$\mathcal{B}$};

\draw [ ->] (10.6, 4.4) --(14.6 ,4.4); 
\draw [ ->] (9, 4) -- (9,2.4);
\draw [ ->] (15.1, 4) -- (15.1,2.4); 
\draw [ ->] (9.5, 2) --(14.7,2); 
 
\end{tikzpicture}
\end{center} 
where $\widehat M$ as the total space of an $S^1$-bundle over the contractible $\mathbb H^2$ can be identified with $\mathbb H^2 \times S^1$. Hence, we can identify the universal cover of $M$ with $\mathbb H^2 \times \mathbb R$, where the vertical lines $\{x\} \times \mathbb R$ universally cover the Seifert fibres of $M$. Denote the normal bundle to these vertical lines by $\tilde \nu$.

The fundamental group of $M$ acts on $\mathbb H^2$ via the composition $\varphi = \bar \varphi \circ \alpha$, which gives rise to a diagonal action $\psi$ of $\pi_1(M)$ on $\widetilde M \times \mathbb H^2$:
$$\psi(\gamma)(x, y) = (\gamma \cdot x, \varphi(\gamma)(M))$$
This action is equivariant with respect to the projection $\widetilde M \times \mathbb H^2 \to \widetilde M$, so taking quotients
determines a locally-trivial $\mathbb H^2$-bundle
$$M \times_\psi \mathbb H^2 = (\widetilde M \times \mathbb H^2)/\pi_1(M) \to M$$

Fix a Riemannian metric on $M$ and pull it back to $\widetilde M$. Since the group of deck transformations of the cover $\widetilde M \to M$ acts as isometries of this metric, the associated
exponential map $\mbox{exp} : T\widetilde M \to \widetilde M$ is $\pi_1(M)$-equivariant, as is the composition 
$$E(\tilde \nu_M) \xrightarrow{\; \mbox{{\footnotesize exp}} \;} \widetilde M = \mathbb H^2 \times  \mathbb R \xrightarrow{\; \mbox{{\footnotesize projection}} \;}  \mathbb H^2$$
Now proceed as in the proof of \cite[Lemma A.3]{BGH25}, starting from its third paragraph, to show that the $\mathbb H^2$-bundle $M \times_\psi \mathbb H^2 \to M$ is topologically equivalent to the normal bundle $\nu_M$ and conclude, as in \cite[Proposition A4]{BGH25}, that $e(\nu_M) = e(\varphi) = e(\bar \varphi \circ \alpha)$. Thus 
\begin{equation}
\label{eqn: Euler normal}
e(\nu_M) = \alpha^*(e(\bar \varphi))
\end{equation}
It is well-known that $e(\bar \varphi)$, respectively $e(\varphi)$, is the obstruction to lifting the representation $\bar \varphi: \Gamma \to PSL_2(\mathbb R)$, respectively $\varphi: \pi_1(M) \to PSL_2(\mathbb R)$, to $\widetilde{PSL}_2(\mathbb R)$ \cite[Lemma 2]{Mil58}. 

An expression for $e(\bar \varphi)$ in terms of the Jankins-Neumann generators $A_0, A_1, \ldots, A_n$ of $H^2(\Gamma)$ can be determined following the method of  \S 3 of \cite{JN85a}. 
To explain it, recall the central extension $1 \to \mathbb Z \to \widetilde{PSL}_2(\mathbb R)  \to PSL_2(\mathbb R) \to 1$. Let $h \in \widetilde{PSL}_2(\mathbb R) $ correspond to $1 \in \mathbb Z$. 

By construction, the generator $\bar x_i$ of the presentation (\ref{eqn: pres 1}) of $\Gamma$ is mapped to a conjugate of a rotation of angle $\frac{2\pi}{b_i}$ in $PSL_2(\mathbb R) \leq \mbox{Homeo}_+(\mathbb S^1)$ by $\bar \varphi$, so lifts to a conjugate of translation by $\frac{1}{b_i}$ in $\widetilde{PSL}_2(\mathbb R)  \leq \mbox{Homeo}_+(\mathbb R)$. Hence, if 
$$1 \to \mathbb Z \to \widetilde \Gamma \to \Gamma \to 1$$
is the pull-back central extension determined by $\bar \varphi$, then $\widetilde \Gamma$ has a presentation 
$$\langle \tilde x_1, \ldots, \tilde x_n,  \tilde w_1, \ldots, \tilde w_g ,  \tilde z_1, \ldots, \tilde z_g  \; | \; \tilde x_i^{b_i} = h, [\tilde w_1, \tilde z_1] \cdots [w_g, \tilde z_g]\tilde x_1 \cdots \tilde x_n = h^{-b} \rangle$$
It is a classic result in the theory of Fuchsian groups that under our assumptions, $b = \chi(|\mathcal{B}|) - n$. As such, the discussion in \S 3 of \cite{JN85a} implies that 
\begin{equation}
\label{eqn: ebarphi}
e(\bar \varphi) = (n - \chi(|\mathcal{B}|))A_0 - A_1 - A_2 - \cdots - A_n
\end{equation}

\begin{proof}[Proof of Theorem \ref{thm: Euler class zero} when $\mathcal{B}$ is hyperbolic]
As $\nu_M$ is an oriented plane bundle, it is trivial if and only if $e(\nu_M) = 0$. Combining (\ref{eqn: Euler normal}), (\ref{eqn: kernel alpha*}), (\ref{eqn: e(E_M)}), and (\ref{eqn: ebarphi}) shows that $e(\nu_M) = 0$ if and only if there is an integer $m$ such that 
\begin{equation}
    \label{equ: hyperbolic iff}
(\chi(|\mathcal{B}|) - n)A_0 + A_1 + \cdots + A_n = m(bA_0 + a_1A_1 + a_2 A_2 + \cdots + a_n A_n)
\end{equation}
Consideration of the normal form for elements of $H^2(\Gamma)$ described in \S \ref{subsec: H^2} shows that the existence of such an $m$ is equivalent to the existence of an integer $m$ for which the following two conditions hold
\vspace{-.2cm}
\begin{enumerate}

\item $ma_i \equiv 1 \mbox{ (mod $b_i$) for } 1 \leq i \leq n$;

\vspace{.2cm} \item $ mb + \sum_{i=1}^n \left\lfloor \frac{ma_i}{b_i} \right\rfloor = \chi(|\mathcal{B}|) - n$. 
\end{enumerate}
For each $i \in \{1, 2, \ldots, n\}$, it follows from the congruence equation $ma_i \equiv 1 \pmod{b_i}$
that 
\[
\left\lfloor \frac{ma_i}{b_i} \right\rfloor 
= \frac{ma_i}{b_i} - \frac{1}{b_i}
\]
Hence, the second equation becomes 
\[
m \left(b + \sum_{i=1}^n \frac{a_i}{b_i}\right) 
= \chi(|\mathcal{B}|) - n + \sum_{i=1}^n \frac{1}{b_i}
\]
which again simplifies to 
\[
m e(M) = \chi(\mathcal{B})
\]
This completes the proof.
\end{proof}

Though we proved Corollary \ref{cor: necessary condition} as a consequence of Theorem \ref{thm: Euler class zero} in \S \ref{subsec: condns}, we provide a proof of the hyperbolic case of the corollary next which is independent of that theorem. 

\begin{proof}[Proof of Corollary \ref{cor: necessary condition} when $\mathcal{B}$ is hyperbolic]
Let $l$ denote the least common multiple of the integers $b_1, b_2, \ldots , b_n$ and note that as $\gcd(\frac{l}{b_1}, \frac{l}{b_2}, \ldots , \frac{l}{b_n}) = 1$, we have an epimorphism 
\begin{equation} 
\label{eqn: epi} 
\theta: H^2(\Gamma) \to \mathbb Z, \; A_0 \mapsto l, A_i \mapsto \frac{l}{b_i} \mbox{ for } 1 \leq i \leq n \nonumber
\end{equation}

Then by (\ref{eqn: e(E_M)}),
\begin{equation}
\label{eqn: Euler number M}
\theta(e(E_M)) = -l \left(b + \sum_{i=1}^n \frac{a_i}{b_i}\right) = -le(M),
\end{equation}
where $e(M)$ is the Euler number of $M$ (see \S \ref{subsection: notation}). 

Next note that (\ref{eqn: ebarphi}) implies that
\begin{equation}
\label{eqn: theta phibar}
\theta(e(\bar \varphi)) = l\left(n - \chi(|\mathcal{B}|) -  \sum_i \frac{1}{b_i}\right) = -l\chi(\mathcal{B})
\end{equation}
Suppose that $e(\nu_M) = 0$. Then (\ref{eqn: Euler normal}) and (\ref{eqn: kernel alpha*}) imply that there is an integer $m$ such that $e(\bar \varphi) = me(E_M) \in H^2(\Gamma)$ and therefore $\theta(e(\bar \varphi)) = m \theta(e(E_M))$. Since $\mathcal{B}$ is hyperbolic, $\chi(\mathcal{B}) \ne 0$, and therefore (\ref{eqn: Euler number M}) and (\ref{eqn: theta phibar}) imply that $m \ne 0$ and $\chi(\mathcal{B}) = me(M)$, which is what was to be proved.  
\end{proof}

\begin{proof}[Proof of Theorem \ref{thm: e = 0 means ctf} when $\mathcal{B}$ is hyperbolic]
We saw in (\ref{eqn: Euler normal}) that that $e(\nu_M) = \alpha^*(e(\bar \varphi)) = e(\bar \varphi \circ \alpha) = e(\varphi)$, where $\varphi = \bar \varphi \circ \alpha: \pi_1(M) \to PSL_2(\mathbb R)$. Hence, if $e(\nu_M) = 0$, then $\varphi$ lifts to a homomorphism $\pi_1(M) \to \widetilde{PSL}_2(\mathbb R) \leq \mbox{Homeo}_+(\mathbb R)$ \cite[Lemma 2]{Mil58}. It follows that $\pi_1(M)$ is left-orderable \cite[Theorem 1.1]{BRW05} and therefore, $M$ admits a horizontal foliation \cite[Theorem 1.3]{BRW05}.
\end{proof}

\section{Seifert fibred examples}
\label{sec: horizontal foliation}
Recall that if $M$ is a Seifert fibred rational homology $3$-sphere, then the underlying surface of the base orbifold is either $S^2$ or $P^2$. However, in the latter case,  $M$ does not admit a co-orientable taut foliation. (See \cite{BRW05} for instance.)

In this section we use Theorem \ref{thm: Euler class zero} and Corollary \ref{cor: necessary condition} to analyse when an oriented Seifert manifold  $M = M(0; b, \frac{a_1}{b_1}, \ldots, \frac{a_n}{b_n})$ has a trivial normal bundle and apply this to show that generically, if $M$ has a co-oriented taut foliation, it has none with zero Euler class. 

We assume that $n \geq 3$. Then the base orbifold $\mathcal{B}$ of $M$ is spherical if $\chi(\mathcal{B}) > 0$, Euclidean if $\chi(\mathcal{B}) = 0$, and hyperbolic if $\chi(\mathcal{B}) < 0$. Orient $\mathcal{B}$ and endow $\nu_M$, the normal bundle to the Seifert fibres of $M$, with the induced orientation. 

\begin{proposition}
\label{prop:ctf if n large}
    Let $M$ be a Seifert fibred rational homology sphere with $n$ exceptional fibres. If $n\geq |H_1(M)|+3$, then $M$ admits a coorientable taut foliation.
\end{proposition}

\begin{proof}
     Suppose that $M$ can be expressed in the form
    \[ M=M\left(0;b, \frac{a_1}{b_1}, \dots, \frac{a_n}{b_n}\right),\]
    where $0<\frac{a_1}{b_1}<1$ for all $i$. We give two proofs, one from the viewpoint of Heegard Floer homology, the other from the viewpoint of foliations transverse to the Seifert fibres of $M$.
    
    {\bf The Heegard Floer approach.} Up to reversing orientation, we can assume that $e(M)<0$. 
    
    Suppose that $M$ does not admit a coorientable taut foliation, or equivalently, that $M$ is an $L$-space  (Lemma \ref{lem: same Euler class} and \cite{LS07}). N\'emethi \cite[\S~11]{Nemethi} provides an algorithm for calculating the Heegaard Floer homology of $M$ in terms of its Seifert invariants. This algorithm shows that if $M$ is an $L$-space then the function
    \[
    \Delta(i)=1-ib - \sum_{j=1}^n \left\lceil \frac{ia_j}{b_j} \right\rceil
    \]
    satisfies $\Delta(i)\geq 0$ for all $i\geq 0$.\footnote{Actually this is equivalent to being an $L$-space.} In particular, one calculates that
    \[
    \Delta(1)=1-b-n\geq 0
    \]
    and that
    \[
    \Delta((b_1\cdots b_n)-1)=1+|H_1(M)|+b\geq 0.
    \]
    Summing these two implies that $n\leq |H_1(M)|+2$, as required.

     {\bf The transverse foliation approach.} After reversing orientation we have $e(M) > 0$, and after reindexing we can suppose that $b_1 \geq b_i$ for each $i$. Then with $p = b_1b_2 \cdots b_n$ and $\sigma = \sum_{i=1}^n a_i/b_i$, we have $p(b + \sigma) = pe(Y) = |H_1(M)|$. Therefore, as $|H_1(M)| \leq n - 3 < nb_2b_3 \cdots b_n$, 
    $$\frac{|H_1(M)|}{p} - n < \frac{|H_1(M)|}{p} - \sigma = b \leq \frac{|H_1(M)|}{p} - \frac{n}{b_1} = \frac{1}{b_1}\left(\frac{|H_1(M)| - nb_2b_3\cdots b_n}{b_2b_3\cdots b_n}\right) < 0. $$
    It follows that $1 - n  \leq b \leq -1$.

    If $2 - n \leq b \leq -2$, $M$ admits a co-orientable taut foliation by \cite[Theorem 2]{JN85b}. 

    If $b = -1$, set $m = p + |H_1(M)| - (n-2)$. Then $2 \leq  m < p$ and if $p_i = p/b_i \in \mathbb Z$ we have  
\begin{itemize}

\item  $a_1p_1 + a_2p_2 + \cdots + a_np_n = p\sigma = p + |H_1(M)| = m + n - 2$; 

\vspace{.2cm} \item  $a_i/b_i = a_ip_i/p < a_ip_i/m$.

\end{itemize}
It follows that the $n$-tuple $(a_1/b_1, \ldots , a_n/b_n)$ lies in $R_n(m)$, where
\begin{eqnarray} 
R_n(m) = \{(\gamma_1, \gamma_2, \ldots, \gamma_n) \in (0,1)^n \;|\; \exists \; (r_1, r_2, \ldots, r_n) \in \mathbb Z^n \mbox{ such that }\nonumber \\
\sum_i r_i = m +n -2 \mbox{ and } \gamma_i < r_i/m \mbox{ for each } i\} \nonumber
\end{eqnarray}
Corollary 3.3 of \cite{JN85b} now implies that $M$ admits a co-orientable taut foliation.

A similar argument can be applied when $b = 1 - n$, after observing that
\begin{eqnarray} 
M = M\left(0;1-n, \frac{a_1}{b_1}, \dots, \frac{a_n}{b_n}\right) & \cong & M\left(0;n-1, -\frac{a_1}{b_1}, \dots, -\frac{a_n}{b_n}\right) \nonumber \\ 
& = & M\left(0; -1, \frac{b_1 - a_1}{b_1}, \ldots, \frac{b_n - a_n}{b_n}\right) \nonumber
\end{eqnarray}
\end{proof}

The following lemma is an immediate consequence of Theorem \ref{thm: Euler class zero}. 

\begin{lemma} 
\label{lemma: gcd}
If $e(\nu_M) = 0$, then $d = \gcd(b_i, b_j)$ implies $a_i \equiv a_j \mbox{ {\rm (mod $d$)}}$. 
\end{lemma} 

\begin{proof}[Proof of Theorem~\ref{thm:all_homology_groups}]
Let $G$ be a finite abelian group of containing an element of order at least $3$. We may write $G\cong \frac{\mathbb{Z}}{d_1}\oplus \dots \oplus \frac{\mathbb{Z}}{d_k}$, where $d_{i+1}$ divides $d_i$ for each $i$ and $d_1\geq 3$. Consider the Seifert fibred space
    \begin{equation}\label{eq:M}
        M=M\left(0;0, \frac{a}{d_1^2}, \frac{1}{d_1}, \dots, \frac{1}{d_k}\right),
    \end{equation}
    where $a=- 1-\sum_{i=1}^k \frac{d_1^2}{d_i}$. For such an $a$ we have $a\equiv -1 \; (\bmod \;d_1)$. Consequently $a$ is coprime to $d_1^2$. Moreover, as $d_1>2$, we see from Lemma \ref{lemma: gcd} that $e(\nu_M) \ne 0$. Calculating the homology using (\ref{eqn: pres 2}) shows that $H_1(M)\cong G$ (compare  \cite[\S~4]{IM20}).

    Next observe that we may perform surgery on the regular fibres in order to obtain a new Seifert fibered space with the same first homology group.
    
    More precisely, suppose that $Y$ is the Seifert fibred rational homology $3$-sphere
    \[ Y=M\left(0;b, \frac{a_1}{b_1}, \dots, \frac{a_n}{b_n}\right)\]
    oriented so that $e(Y) > 0$ and let $H=|H_1(Y)|$ denote the order of the first homology of $Y$. Then one can verify that the manifold
    \[ Y'=M\left(0;b, \frac{a_1}{b_1}, \dots, \frac{a_n}{b_n}, \frac{-H^2}{(b_1\cdots b_n)H +1}\right),\]
    satisfies $H_1(Y)\cong H_1(Y')$. First, one can calculate that 
    \[|H_1(Y')| = (b_1b_2\cdots b_n)(H(b_1b_2\cdots b_n) + 1)e(Y') = |H_1(Y)|.\] On the other hand, consideration of the presentations of $H_1(Y)$ and $H_1(Y')$ obtained from (\ref{eqn: pres 2}) shows that there is an epimorphism $H_1(Y') \to H_1(Y)$ induced by sending $x_i$ to $x_i$ for $1 \leq i \leq n$, $x_{n+1}$ to $1$, and $h$ to $h$. Thus $H_1(Y') \cong H_1(Y)$.  
    
    Now starting with the manifold $M$ as defined in \eqref{eq:M}, and applying the above construction, we obtain an infinite sequence of Seifert fibred spaces, none of which admit a co-orientable taut foliation with Euler class zero by condition (1) of Theorem~\ref{thm: Euler class zero}. It follows from Proposition \ref{prop:ctf if n large} that if 
    \[ Y=M\left(0;b, \frac{a_1}{b_1}, \dots, \frac{a_n}{b_n}\right)\]
    is a Seifert fibred rational homology sphere which is an $L$-space or, equivalently, does not admit a coorientable taut foliation, then $n\leq 2 + |H_1(Y)|$. Thus, the manifolds in this sequence will eventually admit co-orientable taut foliations. 
\end{proof}

Two sets of Seifert invariants $(b, \frac{a_1}{b_1}, \ldots, \frac{a_n}{b_n})$, $(b', \frac{a_1'}{b_1}, \ldots, \frac{a_n'}{b_n})$ are said to be {\it equivalent}, written 
$$\left(b, \frac{a_1}{b_1}, \ldots, \frac{a_n}{b_n}\right) \sim \left(b', \frac{a_1'}{b_1}, \ldots, \frac{a_n'}{b_n}\right)$$ 
if they differ by a finite sequence of moves each of which either multiplies the coordinates by $-1$ or is of the form described in Remark \ref{remark: Seifert notation}(1). Equivalence classes of Seifert invariants classify Seifert structures on a closed $3$-manifold up to homeomorphism. See Remark \ref{remark: Seifert notation}.

\subsection{$\mathcal{B}$ is spherical}
\label{subsec: spherical}
When $\mathcal{B}$ is spherical, $n = 3$ and $(b_1, b_2, b_3)$ is one of the Platonic triples $(2, 2, b_3)$, $(2, 3, 3)$, $(2, 3, 4)$, or $(2, 3, 5)$. In each case, $b_1 = 2$ and $b_2$ is either $2$ or $3$. The fundamental group of $M$ is finite in all cases and therefore $b_1b_2b_3 e(M) = \pm |H_1(M)| \neq 0$ (i.e., $H_1(M)$ is finite). Further, $\chi(\mathcal{B}) > 0$. 

By assumption, $a_i = 1$ when $b_i = 2$. If $b_2 = 3$ and $a_2 = 2$, then the Seifert invariants of $M$ are 
\begin{eqnarray} 
\left(b, \frac{1}{2}, \frac{2}{3}, \frac{a_3}{b_3}\right) \sim \left(b + 1, \frac{1}{2}, \frac{-1}{3}, \frac{a_3}{b_3}\right)  \sim  \left(-b - 1, -\frac{1}{2}, \frac{1}{3}, -\frac{a_3}{b_3}\right) 
\sim  \left(-b - 3, \frac{1}{2}, \frac{1}{3}, \frac{b_3 - a_3}{b_3}\right) \nonumber   
\end{eqnarray}
Thus we can assume that $M$ has Seifert invariants $(b, \frac{1}{2}, \frac{1}{b_2}, \frac{a}{b_3})$. 
The following table lists $\chi(\mathcal{B}), e(M)$ and $|H_1(M)|$ for these manifolds. The last column determines the Seifert invariants of $M$ when $e(\nu_M) = 0$ by applying Theorem \ref{thm: Euler class zero} and Corollary \ref{cor: necessary condition}.

\begin{small} 
\renewcommand{\arraystretch}{1.3} % increase row height (default=1.0)
\begin{center}
\begin{tabular}{|c||c|c|c|c|} \hline 
Seifert invariants &  $\chi(\mathcal{B})$ & $e(M)$ & $\pm|H_1(M)|$ & Seifert invariants when $e(\nu_M) = 0$\\  \hline \hline 
$(b, \tfrac{1}{2}, \tfrac{1}{2}, \tfrac{1}{2})$ & $ \tfrac{1}{2}$ & $\tfrac{2b + 3}{2}$ & $4(2b + 3)$ & $(-1, \tfrac{1}{2}, \tfrac{1}{2}, \tfrac{1}{2})$ and $(-2, \tfrac{1}{2}, \tfrac{1}{2}, \tfrac{1}{2})$ \\  \hline 
$(b, \tfrac{1}{2}, \tfrac{1}{2}, \tfrac{a}{b_3}),\, b_3 \geq 3$ & $\tfrac{1}{b_3}$ & $\tfrac{b_3(b+1) + a}{b_3}$ & $4(b_3(b+1)+a)$ & $(-1, \tfrac{1}{2}, \tfrac{1}{2}, \tfrac{1}{b_3})$ \\  \hline 
$(b, \tfrac{1}{2}, \tfrac{1}{3}, \tfrac{1}{3})$ & $\tfrac{1}{6}$ & $\tfrac{6b+7}{6}$ & $3(6b+7)$ & $(-1, \tfrac{1}{2}, \tfrac{1}{3}, \tfrac{1}{3})$ \\  \hline 
$(b, \tfrac{1}{2}, \tfrac{1}{3}, \tfrac{a}{4})$ & $\tfrac{1}{12}$ & $\tfrac{12b+10+3a}{12}$ & $2(12b+10+3a)$ & $(-1, \tfrac{1}{2}, \tfrac{1}{3}, \tfrac{1}{4})$ \\  \hline 
$(b, \tfrac{1}{2}, \tfrac{1}{3}, \tfrac{a}{5})$ & $\tfrac{1}{30}$ & $\tfrac{30b+6a+25}{30}$ & $30b+6a+25$ & $(-1, \tfrac{1}{2}, \tfrac{1}{3}, \tfrac{1}{5})$ \\  \hline 
\end{tabular}
\end{center}
\end{small}

It is interesting to note that the values of the Seifert invariants when $e(\nu_M) = 0$ correspond to those which realise the minimal values of $|H_1(M)|$ in the given families.  

\begin{remark}
Let $T(r,s)$ be the $(r,s)$-torus knot, where $r < s$ are relatively prime integers greater than 1, and let $g$ denote its genus. If $p/q \ne rs$ then $M = T(r,s)(p/q)$ has a Seifert fibration induced from that on the exterior of $T(r,s)$, and we can use Theorem \ref{thm: Euler class zero} to determine whether or not the Euler class $e(\nu_M)$ is zero. It turns out that if $(r,s) \ne (2,3)$, then $e(\nu_M) = 0$ if and only if $(2g-1)|q| \equiv 1  (\text{{\rm mod} }p)$. Since $ g > 1$, this implies that $p/q < 2g-1$, and therefore $M$ admits a co-orientable taut foliation. On the other hand, for the trefoil $T(2,3)$ we get that $e(\nu_M) = 0$ if and only if $|q| \equiv 1 (\text{{\rm mod} }p)$ and $p/q < 6$. Thus either $p/q < 1$, in which case $M$ has a co-orientable taut foliation, or $p/q = 1,2,3,4$ or $5$. The latter five manifolds $M$ are, respectively, those with Seifert invariants $(-1,1/2,1/3,1/5), (-1,1/2,1/3,1/4), (-1,1/2,1/3,1/3), (-1,1/2,1/2,1/3)$, and $L(5,1)$. As these have finite fundamental groups, they do not admit co-orientable taut foliations.       
\end{remark}

\subsection{$\mathcal{B}$ is Euclidean} 
\label{subsec euclidean} 
Here, either $n = 3$ and $(b_1, b_2, b_3)$ is one of the triples $(2, 3, 6), (2, 4, 4)$, $(3, 3, 3)$ or $n = 4$ and $(b_1, b_2, b_3, b_4) = (2, 2, 2, 2)$. We can assume that $a_i = 1$ in all cases. This is obvious when $b_i = 2$ and using Lemma \ref{lemma: gcd} and the arguments of \S \ref{subsec: spherical} deal with the case that some $b_i > 2$. 

As $\mathcal{B}$ is Euclidean, its Euler characteristic is zero. Hence if $\nu_M$ is trivial, Corollary \ref{cor: necessary condition} implies that $e(M) = 0$ and we can use this to determine the Seifert invariants of $M$ when $e(\nu_M) = 0$.  

\begin{small} 
\renewcommand{\arraystretch}{1.3} % increase row height (default = 1.0)
\begin{center}
\begin{tabular}{|c||c|c|c|} \hline 
Seifert invariants & $e(M)$ & $\pm |H_1(M)|$ & Seifert invariants when $e(\nu_M) = 0$\\  \hline \hline 
$(b, \tfrac{1}{2}, \tfrac{1}{3}, \tfrac{1}{6})$ & $b + 1$ & $36(b+1)$ & $(-1, \tfrac{1}{2}, \tfrac{1}{3}, \tfrac{1}{6})$ \\  \hline 
$(b, \tfrac{1}{2}, \tfrac{1}{4}, \tfrac{1}{4})$ & $b + 1$ & $32(b+1)$ & $(-1, \tfrac{1}{2}, \tfrac{1}{4}, \tfrac{1}{4})$ \\  \hline 
$(b, \tfrac{1}{3}, \tfrac{1}{3}, \tfrac{1}{3})$ & $b + 1$ & $27(b+1)$ & $(-1, \tfrac{1}{3}, \tfrac{1}{3}, \tfrac{1}{3})$ \\  \hline 
$(b, \tfrac{1}{2}, \tfrac{1}{2}, \tfrac{1}{2}, \tfrac{1}{2})$ & $b + 2$ & $16(b+2)$ & $(-2, \tfrac{1}{2}, \tfrac{1}{2}, \tfrac{1}{2}, \tfrac{1}{2})$ \\  \hline 
\end{tabular}
\end{center}
\end{small}

As in the spherical case, those $M$ with $e(\nu_M) = 0$ are characterised by their first homology groups having minimal order ($0$ in this case). Their first Betti numbers are positive and therefore they admit horizontal fibrations (by tori) over the circle. This proves Theorem \ref{thm: e = 0 means ctf} when $\mathcal{B}$ is Euclidean.

\subsection{$\mathcal{B}$ is hyperbolic} 
\label{subsec: hyperbolic}
Throughout this subsection, we assume that $\mathcal B = S^2(b_1,...,b_n)$ is hyperbolic.
We also assume that $M = M(0; b, \frac{a_1}{b_1}, \ldots , \frac{a_n}{b_n})$ is a Seifert manifold with base orbifold $\mathcal{B}$ whose Seifert invariants have been normalised so that $0 < a_i < b_i$ for $1 \le i \le n$.

\begin{remark} 
\label{remark: HF}
The results of \cite{JN85b} and \cite{Na94} show that $M$ admits a horizontal foliation if and only if one of the following conditions holds:
\begin{enumerate}
\setlength{\itemsep}{0.5em}
\item $2 - n \leq b \leq -2$;
\item $b = -1$ and there are coprime integers $0 < c < d$ such that for some permutation $(\frac{c_1}{d}, \frac{c_2}{d}, \ldots, \frac{c_n}{d})$ of $(\frac{c}{d}, \frac{d-c}{d}, \frac{1}{d}, \ldots, \frac{1}{d})$ we have $\frac{a_i}{b_i} < \frac{c_i}{d}$ for each $i$. 
\item $b = 1-n$ and there are coprime integers $0 < c < d$ such that for some permutation $(\frac{c_1}{d}, \frac{c_2}{d}, \ldots, \frac{c_n}{d})$ of $(\frac{c}{d}, \frac{d-c}{d}, \frac{1}{d}, \ldots, \frac{1}{d})$ we have $\frac{b_i - a_i}{b_i} < \frac{c_i}{d}$ for each $i$. 
\end{enumerate}
\end{remark}

\begin{theorem}
\label{thm: hyp base}
 Let $\mathcal B$ be as above.
 \begin{enumerate}[leftmargin=*]
 \setlength\itemsep{0.5em}
     \item[\rm (1)] There exists a Seifert fibred rational homology sphere $M$ with base orbifold $\mathcal B$ such that $M$ has a horizontal foliation and $e(\nu_M) = 0$.
\item[\rm (2)] With exactly 58 exceptions $\mathcal B$, there exists a Seifert fibred rational homology sphere $M$ with base orbifold $\mathcal B$ such that $M$ has a horizontal foliation and $e(\nu_{M}) \ne 0$.
 \end{enumerate}
\end{theorem}

Thus the manifolds $M$ in part (2) have co-orientable taut foliations but none with Euler class zero. 

To prove part (1) we take $M = M(0;2-n,1/b_1,...,1/b_n)$. It follows from Remark \ref{remark: HF}  that $M$ has a horizontal foliation. (When $n=3$, use the fact that $\mathcal B$ is hyperbolic.) Taking $m=1$ in Theorem \ref{thm: Euler class zero} shows that $e(\nu_M) = 0$.

The exceptions in part (2) are explicitly listed in the Appendix. We will not include a complete proof of part (2); we content ourselves with proving two partial results, Theorems \ref{thm: large n} and \ref{thm: special n = 3} below, to illustrate the general strategy.

First we introduce some notation. Let $X = -\chi(\mathcal B) = (n-2) - \sum_{i = 1}^n {1/b_i}$. Since $\mathcal{B}$ is hyperbolic, $X > 0$. Next let $M$ be a Seifert fibre space with base orbifold $\mathcal B$; thus $M$ is of the form $M(0;b,a_1/b_1,...,a_n/b_n)$. Then
$$-e(M) = -b - \sum_{i = 1}^n {a_i/b_i} = X - \left((n-2) + b + \sum_{i = 1}^n \frac{a_i-1}{b_i}\right)$$
Let $\alpha = (n-2) + b + \sum_{i = 1}^n {(a_i-1)/b_i}$. Then $$-e(M) = X -\alpha$$

By Theorem \ref{thm: Euler class zero}, $e(\nu_M) = 0$ if and only if there exists an integer $m$ such that $m a_i \equiv 1$ (mod $b_i$) for $1 \le i \le n$ and $X = m(X-\alpha)$, i.e. $m \alpha = (m-1)X$. We note that if $X-\alpha > 0$ then $M$ is a rational homology sphere, and if in addition $\alpha \ne 0$ then $m > 1$.

\begin{theorem}
\label{thm: large n}  
  Let $\mathcal B = S^2(b_1,...,b_n)$ be as above, where $n \ge 7$. If $n = 7$, assume that some $b_i$ is greater than $2$. Then there is a Seifert fibred rational homology sphere $M$ with base orbifold $\mathcal B$ such that $M$ has a horizontal foliation and $e(\nu_M) \ne 0$.   
\end{theorem}

\begin{proof}
  We take $M$ to be $M(0;3-n,1/b_1,...,1/b_n)$. Then $M$ has a horizontal foliation by Remark \ref{remark: HF}.

  We have $X \ge (n-2) - n/2 = (n-4)/2$, with strict inequality if $n = 7$. Hence $X > 3/2$. Also, $\alpha = 1$, so $M$ is a rational homology sphere.

  If $e(\nu_M) = 0$ then there exists an integer $m > 1$ such that $m = m \alpha = (m-1)X > (3/2)(m-1)$. Therefore $m < 3$, so $m = 2$. But then $m \not\equiv 1$ (mod $b_i$) for any $i$.
\end{proof}

The exclusion of the case  $\mathcal B = S^2(2,2,...,2)$, $n = 7$, is necessary. To see this, note that here $X = (7-2) - 7/2 = 3/2$. The manifolds $M$ with base orbifold $\mathcal B$ are those of the form $M(0;b,1/2,1/2,...,1/2)$. By Remark \ref{remark: HF} $M$ has a horizontal foliation if and only if $-5 \le b \le -2$, i.e. $\alpha = k$, $0 \le k \le 3$. The corresponding manifold $M_k$, say, has $e(\nu_{M_k}) = 0$ if and only if there exists an odd integer $m_k$ such that $X = m_k(X-k)$. Hence $m_k = 3/(3-2k)$, so $m_0 = 1, m_1 = 3, m_ 2 = -3$, and $m_3 = -1.$ (In fact $M_{3-k} \cong - M_k.)$ Thus every $M$ with base orbifold $\mathcal B$ that has a horizontal foliation has $e(\nu_M) = 0$, so $\mathcal B(2,2,...,2)$, $n = 7$, is one of the exceptions in part (2) of Theorem \ref{thm: hyp base}. 

The next theorem deals with the case $n=3$ when each $b_i$ is at least $5$.

\begin{theorem}
\label{thm: special n = 3}
  Let $\mathcal B = S^2(b_1,b_2,b_3)$, where $b_i \ge 5$, $1 \le i \le 3$, and the $b_i$ are not all $6$. Then there is a Seifert fibred rational homology sphere $M$ with base orbifold $\mathcal B$ such that $M$ has a horizontal foliation and $e(\nu_M) \ne 0$.
\end{theorem}

The proof will use the following lemma.

\begin{lemma}
\label{lemma: BP}
If either $c = 5$ or $c \ge 7$ then there exists a prime $p$ which does not divide $c$ such that $(p-1)/c \le 1/3$ and $p/c < 1/2$.
\end{lemma}

\begin{proof}
  If $c$ is odd take $p = 2$.

  If $c$ is even and not a multiple of $3$, take $p = 3$. 

  So we may suppose that $c = 6n$, where $n > 1$. By Bertrand's Postulate there exists a prime $p$ such that $n+2 < p < 2n+2$. The left-hand inequality implies that $p$ does not divide $c$ and the right-hand inequality shows that $p-1 \le 2n$. Therefore $(p-1)/c \le 2n/6n = 1/3$. Also $p/c \le 1/3 + 1/c < 1/2$. 
\end{proof}

\begin{proof}[Proof of Theorem \ref{thm: special n = 3}]
  If two of the $b_i$, say $b_1$ and $b_2$, are equal to $5$ then we take $M = M(0;-1,1/5,2/5,1/b_3)$; see Lemma \ref{lemma: gcd}. We may therefore assume that $X > 1 - 3/5 = 2/5.$

  By hypothesis, some $b_i$, say $b_1$, is not equal to $6$. Then by Lemma \ref{lemma: BP} there is a prime $p$ such that $(p-1)/b_1 \le 1/3$ and $p/b_1 < 1/2$. Let $M = M(0;-1,p/b_1,1/b_2,1/b_3)$. Since $p/b_1 < 1/2$, $M$ has a horizontal foliation by Remark \ref{remark: HF}. Also, $\alpha = (p-1)/b_1$, so $X - \alpha > 0$. Hence $M$ is a rational homology sphere. If $e(\nu_M) = 0$ then there exists $m > 1$ such that $m\alpha = (m-1)X$. Therefore $m/3 > (2/5)(m-1)$. Hence $5m > 6(m-1)$, so $m < 6$. But then, since $b_2 \ge 5$, $m \not\equiv 1$ (mod $b_2$).
\end{proof}

It is necessary to exclude the case $\mathcal B = S^2(6,6,6)$. For it is easy to see from Remark \ref{remark: HF} that the only manifold $M$ with base orbifold $\mathcal B$ that has a horizontal foliation is (up to orientation) $M(0;-1,1/6,1/6,1/6)$, which has $e(\nu_M) = 0$. 

\bigskip

{
\footnotesize
\bibliographystyle{alpha}
 %plain
\bibliography{bgh_zero}
}

\newpage
\appendix
\section{Exceptions to Theorem \ref {thm: hyp base}}

Here we list the orbifolds $S^2(b_1,...,b_n)$ that are exceptions in part (2) of Theorem \ref {thm: hyp base}. We adopt the notational convention stated at the beginning of \S \ref{subsec: hyperbolic}, and  also assume that $b_1 \le \cdots \le b_n$.

Let $M = M(0;b,a_1/b_1,...,a_n/b_n)$ be a Seifert fibre space with base orbifold $S^2(b_1,...,b_n)$. We distinguish three types of exceptional $S^2(b_1,...,b_n)$.

\begin{enumerate}[leftmargin=*] 
\setlength\itemsep{0.5em}
    \item $M$ has a horizontal foliation if and only if up to orientation, i.e. for either $M$ or $-M$, $a_1 = \ldots = a_n = 1$ and $b = 2-n$. Then taking $m=1$ shows that $e(\nu_M) = 0$.
    \item  There exists $M$ with a horizontal foliation such that for both $M$ and $-M$ either $b \ne 2-n$ or some $a_i \ne 1$, but any such $M$ has $e(\nu_M) = 0$.
    \item There exists $M$ with a horizontal foliation having $e(\nu_M) \ne 0$, but any such $M$ has $e(M) = 0$, and hence is not a rational homology sphere. 
\end{enumerate}

The exceptions of Type (1) can be easily determined using Remark \ref{remark: HF}. They are  
\begin{itemize}[leftmargin=*]
 \setlength\itemsep{0.5em}
    \item $S^2(2,3,c)$, where $c = 7, 8, 9, 10, 12, 14, 18, 24,$ or $30$
    \item $S^2(2,4,c)$, where $c = 5, 6, 8,$ or $12$
    \item $S^2(2,5,c)$, where $c = 5$ or $6$
    \item $S^2(2,6,c)$, where $c = 6$ or $12$
    \item $S^2(3,3,c)$, where $c = 4$ or $6$
    \item $S^2(3,4,c)$, where $c = 4$ or $6$
    \item $S^2(4,4,c)$, where $c = 4$ or $6$
    \item $S^2(6,6,6)$
    \item $S^2(2,2,2,c)$, where $c = 3, 4,$ or $6$
    \item $S^2(2,2,2,2,2)$
\end{itemize}

For the exceptions of Type (2), we list the corresponding manifolds $M$ up to orientation, and the values of $m$ which show that $e(\nu_M) = 0$.
\begin{itemize}[leftmargin=*, itemsep=5pt]
    \item $S^2(2,3,c)$, where $c = 11,13,16$, or $20$ and $M = M(0;-1,1/2,1/3,a/c)$, $1 < a < c/2$. In particular, 
    \begin{itemize}
    \item when $c = 11$, $a = 2$ and $m = -5$
    \item when $c = 13$, $a = 2$ and $m = 7$
    \item when $c = 16$, $a = 3$ and $m = -5$
    \item when $c = 20$, $a = 3$ and $m = 7$
\end{itemize}
\item $S^2(2,4,c)$, where $c = 7, 9, 10, 14, 18,$ or $30$ and $M = M(0;-1,1/2,1/4,a/c)$, $1 < a < c/2$. In particular, 
\begin{itemize}
    \item when $c = 7$, $a = 2$ and $m = -3$
    \item when $c = 9$, $a = 2$ and $m = 5$
    \item when $c = 10$, $a = 3$ and $m = -3$
    \item when $c = 14$, $a = 3$ and $m = 5$
    \item when $c = 18$, $a = 3$ and $m = -7$
    \item when $c = 30$, $a = 7$ and $m = 13$
\end{itemize}
\item $S^2(2,6,c)$, where $c = 8$ or $10$ and $M = M(0;-1,1/2,1/6,a/c)$, $1 < a < c/2$. In particular, 
\begin{itemize}
    \item when $c = 8$, $a = 3$ and $m = -5$
    \item when $c = 10$, $a = 3$ and $m = 7$
\end{itemize}
\item $S^2(3,3,c)$, where $c = 5, 7, 8,$ or $10$ and $M = M(0;-1,1/3,1/3,a/c)$, $1 < a < c/2$. In particular, 
\begin{itemize}
    \item when $c = 5$, $a = 2$ and $m = -2$
    \item when $c = 7$, $a = 2$ and $m = -5$; or  $a = 3$ and $m = -2$
    \item when $c = 8$, $a = 3$ and $m = -5$
    \item when $c = 10$, $a = 3$ and $m = 7$
\end{itemize}
\item $S^2(3,4,5)$, where $M = M(0;-1,1/3,1/4,2/5)$ and $m = 13$
\item $S^2(4,4,12)$, where $M = M(0;-1,1/4,1/4,5/12)$ and $m = 5$

\item $S^2(4,6,6)$, where $M = M(3/4,1/6,1/6)$ and $m = -5$
\item $S^2(2,2,2,c)$, where $c = 5, 8,$ or $12$ and $M = M(0;-2,1/2,1/2,1/2,a/c)$, $1 < a < c/2$. In particular, 
\begin{itemize}
    \item when $c = 5$, $a = 2$ and $m = 3$
    \item when $c = 8$, $a = 3$ and $m = 3$
    \item when $c = 12$, $a = 5$ and $m = 5$
\end{itemize}
\item $S^2(2,2,3,4)$, where $M = M(0;-2,1/2,1/2, 2/3,1/4)$ and $m = 5$
\item $S^2(2,2,4,6)$, where $M = M(0;-2,1/2,1/2,3/4,1/6)$ and $m = 7$
\item $S^2(2,2,2,2,c)$, where $c = 4$ or $6$. In particular, 
\begin{itemize}
    \item when $c = 4$, $M = M(0;-2,1/2,1/2,1/2,1/2,1/4)$ and $m = -3$
    \item when $c = 6$, $M = M(0;-2,1/2,1/2,1/2,1/2,1/6)$ and $m = -5$
\end{itemize}
\item $S^2(2,2,2,2,2,2,2)$. See \S \ref{subsec: hyperbolic}.
\end{itemize}

Lastly, we list the exceptions of Type~(3). In each case, we list the only manifolds $M$, up to orientation, that admit a horizontal foliation with $e(\nu_M) \ne 0$. One easily checks that in all cases $e(M) = 0$: 

\begin{itemize}[itemsep=6pt, leftmargin=*]
    \item $S^2(2,5,10)$, where $M = M(0;-1,1/2,1/5,3/10)$ and $M = M(0;-1,1/2,2/5,1/10)$
    \item $S^2(3,4,12)$, where $M = M(0;-1,1/3,1/4,5/12)$
    \item $S^2(2,2,2,2,2,2)$, where $M = M(0;-3,1/2,1/2,1/2,1/2,1/2,1/2)$
\end{itemize}

\newpage

\end{document}